\documentclass[11pt, a4paper]{article}
\usepackage[cp1251]{inputenc}
\usepackage[russian]{babel}
\usepackage{amsmath} \usepackage{euscript}
\usepackage{mymatrx5}
\def\mymatrix{\MyMatrixwithdelims..}
\Linestrue\Autonumtrue
\usepackage{mathrsfs}
\usepackage{amssymb}
\begin{document}

\setcounter{MaxMatrixCols}{14}
\newcounter{bnomer} \newcounter{snomer}
\newcounter{bsnomer}
\setcounter{bnomer}{0}
\renewcommand{\thesnomer}{\thebnomer.\arabic{snomer}}
\renewcommand{\thebsnomer}{\thebnomer.\arabic{bsnomer}}
\renewcommand{\refname}{\begin{center}\large{\textbf{References}}\end{center}}

\newcommand{\sect}[1]{%
\setcounter{snomer}{0}\setcounter{bsnomer}{0}
\refstepcounter{bnomer}
\par\bigskip\begin{center}\large{\textbf{\arabic{bnomer}. {#1}}}\end{center}}
\newcommand{\sst}[1]{%
\refstepcounter{bsnomer}
\par\bigskip\textbf{\arabic{bnomer}.\arabic{bsnomer}. {#1}.}}
\newcommand{\defi}[1]{%
\refstepcounter{snomer}
\par\medskip\textbf{Definition \arabic{bnomer}.\arabic{snomer}. }{#1}\par\medskip}
\newcommand{\theo}[2]{%
\refstepcounter{snomer}
\par\textbf{Теорема \arabic{bnomer}.\arabic{snomer}. }{#2} {\emph{#1}}\hspace{\fill}$\square$\par}
\newcommand{\mtheop}[2]{%
\refstepcounter{snomer}
\par\textbf{Theorem \arabic{bnomer}.\arabic{snomer}. }{\emph{#1}}
\par\textsc{Proof}. {#2}\hspace{\fill}$\square$\par}
\newcommand{\mcorop}[2]{%
\refstepcounter{snomer}
\par\textbf{Corollary \arabic{bnomer}.\arabic{snomer}. }{\emph{#1}}
\par\textsc{Proof}. {#2}\hspace{\fill}$\square$\par}
\newcommand{\mtheo}[1]{%
\refstepcounter{snomer}
\par\medskip\textbf{Theorem \arabic{bnomer}.\arabic{snomer}. }{\emph{#1}}\par\medskip}
\newcommand{\mlemm}[1]{%
\refstepcounter{snomer}
\par\medskip\textbf{Lemma \arabic{bnomer}.\arabic{snomer}. }{\emph{#1}}\par\medskip}
\newcommand{\mprop}[1]{%
\refstepcounter{snomer}
\par\medskip\textbf{Proposition \arabic{bnomer}.\arabic{snomer}. }{\emph{#1}}\par\medskip}
\newcommand{\theobp}[2]{%
\refstepcounter{snomer}
\par\textbf{Теорема \arabic{bnomer}.\arabic{snomer}. }{#2} {\emph{#1}}\par}
\newcommand{\theop}[2]{%
\refstepcounter{snomer}
\par\textbf{Theorem \arabic{bnomer}.\arabic{snomer}. }{\emph{#1}}
\par\textsc{Proof}. #2\hspace{\fill}$\square$\par}
\newcommand{\theosp}[2]{%
\refstepcounter{snomer}
\par\textbf{Теорема \arabic{bnomer}.\arabic{snomer}. }{\emph{#1}}
\par\textbf{Схема доказательства}. {#2}\hspace{\fill}$\square$\par}
\newcommand{\exam}[1]{%
\refstepcounter{snomer}
\par\medskip\textbf{Example \arabic{bnomer}.\arabic{snomer}. }{#1}\par\medskip}
\newcommand{\deno}[1]{%
\refstepcounter{snomer}
\par\textbf{Definition \arabic{bnomer}.\arabic{snomer}. }{#1}\par}
\newcommand{\post}[1]{%
\refstepcounter{snomer}
\par\textbf{Предложение \arabic{bnomer}.\arabic{snomer}. }{\emph{#1}}\hspace{\fill}$\square$\par}
\newcommand{\postp}[2]{%
\refstepcounter{snomer}
\par\medskip\textbf{Proposition \arabic{bnomer}.\arabic{snomer}. }{\emph{#1}}%
\ifhmode\par\fi\textsc{Proof}. {#2}\hspace{\fill}$\square$\par\medskip}
\newcommand{\lemm}[1]{%
\refstepcounter{snomer}
\par\textbf{Lemma \arabic{bnomer}.\arabic{snomer}. }{\emph{#1}}\hspace{\fill}$\square$\par}
\newcommand{\lemmp}[2]{%
\refstepcounter{snomer}
\par\medskip\textbf{Lemma \arabic{bnomer}.\arabic{snomer}. }{\emph{#1}}
\par\textsc{Proof}. #2\hspace{\fill}$\square$\par\medskip}
\newcommand{\coro}[1]{%
\refstepcounter{snomer}
\par\textbf{Corollary \arabic{bnomer}.\arabic{snomer}. }{\emph{#1}}\hspace{\fill}$\square$\par}
\newcommand{\mcoro}[1]{%
\refstepcounter{snomer}
\par\textbf{Corollary \arabic{bnomer}.\arabic{snomer}. }{\emph{#1}}\par\medskip}
\newcommand{\corop}[2]{%
\refstepcounter{snomer}
\par\textbf{Corollary \arabic{bnomer}.\arabic{snomer}. }{\emph{#1}}
\par\textsc{Proof}. {#2}\hspace{\fill}$\square$\par}
\newcommand{\nota}[1]{%
\refstepcounter{snomer}
\par\medskip\textbf{Remark \arabic{bnomer}.\arabic{snomer}. }{#1}\par\medskip}
\newcommand{\propp}[2]{%
\refstepcounter{snomer}
\par\medskip\textbf{Proposition \arabic{bnomer}.\arabic{snomer}. }{\emph{#1}}
\par\textsc{Proof}. {#2}\hspace{\fill}$\square$\par\medskip}
\newcommand{\hypo}[1]{%
\refstepcounter{snomer}
\par\medskip\textbf{Conjecture \arabic{bnomer}.\arabic{snomer}. }{\emph{#1}}\par\medskip}
\newcommand{\prop}[1]{%
\refstepcounter{snomer}
\par\textbf{Proposition \arabic{bnomer}.\arabic{snomer}. }{\emph{#1}}\hspace{\fill}$\square$\par}

\newcommand{\Ind}[3]{%
\mathrm{Ind}_{#1}^{#2}{#3}}
\newcommand{\Res}[3]{%
\mathrm{Res}_{#1}^{#2}{#3}}
\newcommand{\epsi}{\epsilon}
\newcommand{\tri}{\triangleleft}
\newcommand{\Supp}[1]{%
\mathrm{Supp}(#1)}

\makeatletter
\def\iddots{\mathinner{\mkern1mu\raise\p@
\vbox{\kern7\p@\hbox{.}}\mkern2mu
\raise4\p@\hbox{.}\mkern2mu\raise7\p@\hbox{.}\mkern1mu}}
\makeatother

\newcommand{\reg}{\mathrm{reg}}
\newcommand{\empr}[2]{[-{#1},{#1}]\times[-{#2},{#2}]}
\newcommand{\sreg}{\mathrm{sreg}}
\newcommand{\codim}{\mathrm{codim}\,}
\newcommand{\chara}{\mathrm{char}\,}
\newcommand{\rk}{\mathrm{rk}\,}
\newcommand{\chr}{\mathrm{ch}\,}
\newcommand{\id}{\mathrm{id}}
\newcommand{\Ad}{\mathrm{Ad}}
\newcommand{\Ker}{\mathrm{Ker}\,}
\newcommand{\End}{\mathrm{End}}
\newcommand{\Ann}{\mathrm{Ann}\,}
\newcommand{\col}{\mathrm{col}}
\newcommand{\row}{\mathrm{row}}
\newcommand{\low}{\mathrm{low}}
\newcommand{\pho}{\hphantom{\quad}\vphantom{\mid}}
\newcommand{\fho}[1]{\vphantom{\mid}\setbox0\hbox{00}\hbox to \wd0{\hss\ensuremath{#1}\hss}}
\newcommand{\wt}{\widetilde}
\newcommand{\wh}{\widehat}
\newcommand{\ad}[1]{\mathrm{ad}_{#1}}
\newcommand{\tr}{\mathrm{tr}\,}
\newcommand{\GL}{\mathrm{GL}}
\newcommand{\Pf}{\mathrm{Pf}}
\newcommand{\Prim}{\mathrm{Prim}\,}
\newcommand{\Cent}{\mathrm{Cent}\,}
\newcommand{\SL}{\mathrm{SL}}
\newcommand{\SO}{\mathrm{SO}}
\newcommand{\Sp}{\mathrm{Sp}}
\newcommand{\Mat}{\mathrm{Mat}}
\newcommand{\Stab}{\mathrm{Stab}}
\newcommand{\ilm}{\varinjlim}
\newcommand{\lee}{\leq}
\newcommand{\gee}{\geq}

\newcommand{\vfi}{\varphi}
\newcommand{\vpi}{\varpi}
\newcommand{\teta}{\vartheta}
\newcommand{\Bfi}{\Phi}
\newcommand{\Fp}{\mathbb{F}}
\newcommand{\Rp}{\mathbb{R}}
\newcommand{\Zp}{\mathbb{Z}}
\newcommand{\Cp}{\mathbb{C}}
\newcommand{\Np}{\mathbb{N}}
\newcommand{\Xt}{\mathfrak{X}}
\newcommand{\ut}{\mathfrak{u}}
\newcommand{\at}{\mathfrak{a}}
\newcommand{\nt}{\mathfrak{n}}
\newcommand{\mt}{\mathfrak{m}}
\newcommand{\htt}{\mathfrak{h}}
\newcommand{\spt}{\mathfrak{sp}}
\newcommand{\slt}{\mathfrak{sl}}
\newcommand{\ot}{\mathfrak{so}}
\newcommand{\rt}{\mathfrak{r}}
\newcommand{\rad}{\mathfrak{rad}}
\newcommand{\bt}{\mathfrak{b}}
\newcommand{\gt}{\mathfrak{g}}
\newcommand{\vt}{\mathfrak{v}}
\newcommand{\pt}{\mathfrak{p}}
\newcommand{\Po}{\mathcal{P}}
\newcommand{\Uo}{\EuScript{U}}
\newcommand{\Fo}{\EuScript{F}}
\newcommand{\Do}{\EuScript{D}}
\newcommand{\Eo}{\EuScript{E}}
\newcommand{\Iu}{\mathcal{I}}
\newcommand{\Du}{\mathcal{D}}
\newcommand{\Uu}{\mathcal{U}}
\newcommand{\Mo}{\mathcal{M}}
\newcommand{\Nu}{\mathcal{N}}
\newcommand{\Ro}{\mathcal{R}}
\newcommand{\Co}{\mathcal{C}}
\newcommand{\Lo}{\mathcal{L}}
\newcommand{\Ou}{\mathcal{O}}
\newcommand{\Au}{\mathcal{A}}
\newcommand{\Vu}{\mathcal{V}}
\newcommand{\Bu}{\mathcal{B}}
\newcommand{\Sy}{\mathcal{Z}}
\newcommand{\Sb}{\mathcal{F}}
\newcommand{\Gr}{\mathcal{G}}
\newcommand{\rtc}[1]{C_{#1}^{\mathrm{red}}}

\author{Mikhail Ignatyev\and Ivan Penkov}

\date{}
\title{\Large{Infinite Kostant cascades and centrally generated\\ primitive ideals of $U(\nt)$ in types $A_{\infty}$, $C_{\infty}$\mbox{$\vphantom{1}$}\footnotetext{The first author has been supported in part by RFBR grants no. 14--01--31052 and 14--01--97017, by the Dynasty Foundation and by the Ministry of Science and Education of the Russian Federation. A part of this work was done in the Max Planck Institute for Mathematics in Bonn and in Jacobs University Bremen, and the first author thanks these institutions for their hospitality. The second author has been supported in part by DFG via the priority program Representation Theory.}}
}\maketitle

\begin{center}
\begin{tabular}{p{15cm}}
\small{\textsc{Abstract}. We study the center of $U(\nt)$, where $\nt$ is the locally nilpotent radical of a splitting Borel subalgebra of a simple complex Lie algebra $\gt=\slt_{\infty}(\Cp)$, $\ot_{\infty}(\Cp)$, $\spt_{\infty}(\Cp)$. There are infinitely many isomorphism classes of Lie algebras $\nt$, and we provide explicit generators of the center of $U(\nt)$ in all cases. We then fix $\nt$ with ``largest possible'' center of $U(\nt)$ and characterize the centrally generated primitive ideals of $U(\nt)$ for $\gt=\slt_{\infty}(\Cp)$, $\spt_{\infty}(\Cp)$ in terms of the above generators. As a~preliminary result, we provide a characterization of the centrally generated primitive ideals in the enveloping algebra of the nilradical of a Borel subalgebra of $\slt_n(\Cp)$, $\spt_{2n}(\Cp)$.}
\end{tabular}
\end{center}

\sect{Introduction}

The theory of primitive ideals in enveloping algebras of Lie algebras has its roots in the re\-pre\-sen\-tation theory of Lie algebras. However, classifying irreducible representations of Lie algebras is not feasible except in few very special cases, while a classification of annihilators of irreducible representations, i.e., of primitive ideals, can be achieved in much greater generality. This idea goes back to J. Dixmier and his seminar, and for semisimple or solvable finite-dimensional Lie algebras there is an extensive theory of primitive ideals.

In the case when $\nt$ is a finite-dimensional nilpotent Lie algebra, the primitive ideals in the universal enveloping algebra $U(\nt)$ can be described in terms of the Dixmier map assigning to any linear form $f\in\nt^*$ a primitive ideal $J(f)$ of $U(\nt)$. If $\nt$ is abelian, $J(f)$ is simply the annihilator of $f$. For a general finite-dimensional nilpotent Lie algebra $\nt$, the theory of primitive ideals retains many properties from the abelian case: in particular, $J(f)$ is always a maximal ideal and every primitive ideal in $U(\nt)$ is of the form $J(f)$ for some $f\in\nt^*$. Moreover, $J(f)=J(f')$ if and only if $f$ and $f'$ belong to the same coadjoint orbit in $\nt^*$.

The idea of classifying primitive ideals rather then irreducible representations makes even more sense for infinite-dimensional Lie algebras, and in this paper we make some first steps in this direction for a natural class of locally nilpotent infinite-dimensional Lie algebras. These are the locally nilpotent radicals $\nt$ of splitting Borel subalgebras of the three simple finitary complex Lie algebras $\slt_{\infty}(\Cp)$, $\ot_{\infty}(\Cp)$, $\spt_{\infty}(\Cp)$.

A comprehensive theory of primitive ideals in $U(\nt)$ remains to be built. In this paper we concentrate on centrally generated primitive ideals in $U(\nt)$. We first provide a description of the center of the enveloping algebra $U(\nt)$ for any locally nilpotent radical $\nt$ as above, and then use the result to describe the centrally generated primitive ideals in $U(\nt)$ for some interesting choices of $\nt$.

The splitting Borel subalgebras of $\gt=\slt_{\infty}(\Cp)$, $\ot_{\infty}(\Cp)$, $\spt_{\infty}(\Cp)$ are not conjugate, and there are infinitely many isomorphism classes of locally nilpotent radicals $\nt$. In the finite-dimensional case Kostant cascades of orthogonal roots play an important role in describing the center of $U(\nt)$. In the infinite-dimensional case the center of $U(\nt)$ is described in terms of a possibly infinite Kostant cascade. A significant difference with the finite-dimensional case is that the Kostant cascade depends in an essential way on the isomorphism class of $\nt$, and that in most cases the cascade is finite rather than infinite.

In order to obtain an explicit form of the generators of the center of $U(\nt)$, we first recall such an explicit form of the generators in the finite-dimensional case due to A.~Joseph, F.~Fauquant-Millet, R.~Lipsman, J.A.~Wolf, and A.~Panov. This enables us to give an explicit description of the center of $U(\nt)$ in all cases.

We then concentrate on the case when $\gt=\slt_{\infty}(\Cp)$, $\spt_{\infty}(\Cp)$ and $U(\nt)$ has the ``largest possible'' center. This latter requirement singles out only one isomorphism class of subalgebras $\nt$ for a fixed~$\gt$. For such $\nt$ we construct a Dixmier map defined for certain linear forms $f\in\nt^*$, closely related to the Kostant cascade of $\nt$. We refer to these forms as Kostant forms. Our main result implies then that the Dixmier map establishes a one-to-one correspondence between Kostant forms and centrally generated primitive ideals in $U(\nt)$. This provides an explicit description of the centrally generated primitive ideals of $U(\nt)$. As a corollary we obtain that centrally generated primitive ideals $J$ of $U(\nt)$ are maximal ideals, and that the quotient $U(\nt)/J$ is a Weyl algebra with infinitely many generators.

To the best of our knowledge, the analogous description, via Kostant forms, of the centrally generated primitive ideals in $U(\nt)$ for $\nt\subset\slt_n(\Cp)$, $\spt_{2n}(\Cp)$ is also new.

We thank A. Joseph and A. Panov for helpful discussions.

\sect{The center of $U(\nt)$}


\sst{Finite-dimensional\label{sst:centers_finite} case} Let $n\in\Zp_{>0}$. Throughout this subsection $\gt$ denotes one of the Lie algebras $\slt_n(\Cp)$, $\ot_{2n}(\Cp)$, $\ot_{2n+1}(\Cp)$ or $\spt_{2n}(\Cp)$. The algebra $\ot_{2n}(\Cp)$ (respectively, $\ot_{2n+1}(\Cp)$ and $\spt_{2n}(\Cp)$) is realized as the sub\-al\-gebra of $\slt_{2n}(\Cp)$ (respectively, $\slt_{2n+1}(\Cp)$ and $\slt_{2n}(\Cp)$) consisting of all~$x$ such that $\beta(u,xv)+\beta(xu,v)=0$ for all $u,v$ in $\Cp^{2n}$ (respectively, in $\Cp^{2n+1}$ and $\Cp^{2n}$), where
\begin{equation*}
\beta(u,v)=\begin{cases}
\sum\nolimits_{i=1}^n(u_iv_{-i}+u_{-i}v_i)&\text{for }\ot_{2n}(\Cp),\\
u_0v_0+\sum\nolimits_{i=1}^n(u_iv_{-i}+u_{-i}v_i)&\text{for }\ot_{2n+1}(\Cp),\\
\sum\nolimits_{i=1}^n(u_iv_{-i}-u_{-i}v_i)&\text{for }\spt_{2n}(\Cp).
\end{cases}
\end{equation*}
Here for $\ot_{2n}(\Cp)$ (respectively, for $\ot_{2n+1}$ and $\spt_{2n}(\Cp))$ we denote by $e_1,\ldots,e_n,e_{-n},\ldots,e_{-1}$ (res\-pec\-tively, by $e_1,\ldots,e_n,e_0,e_{-n},\ldots,e_{-1}$ and $e_1,\ldots,e_n,e_{-n},\ldots,e_{-1}$) the standard basis of $\Cp^{2n}$ (res\-pec\-tively, of $\Cp^{2n+1}$ and $\Cp^{2n}$), and by $x_i$ the coordinate of a vector $x$ corresponding to $e_i$.

The set of all diagonal matrices from $\gt$ is a Cartan subalgebra of $\gt$; we denote it by $\htt$. Let $\Phi$ be the root system of $\gt$ with respect to $\htt$. Note that $\Phi$ is of type $A_{n-1}$ (respectively, $D_n$, $B_n$ and $C_n$) for $\slt_n(\Cp)$ (respectively, for $\ot_{2n}(\Cp)$, $\ot_{2n+1}(\Cp)$ and $\spt_{2n}(\Cp)$). The set of all upper-triangular matrices from $\gt$ is a Borel subalgebra of $\gt$ containing $\htt$; we denote it by $\bt$. Let $\Phi^+$ be the set of positive roots with respect to $\bt$. As usual, we identify $\Phi^+$ with the following subset of $\Rp^n$:
\begin{equation*}
\begin{split}
A_{n-1}^+&=\{\epsi_i-\epsi_j,~1\leq i<j\leq n\},\\
B_n^+&=\{\epsi_i-\epsi_j,~1\leq i<j\leq n\}\cup\{\epsi_i+\epsi_j,~1\leq i<j\leq n\}\cup\{\epsi_i,~1\leq i\leq n\},\\
C_n^+&=\{\epsi_i-\epsi_j,~1\leq i<j\leq n\}\cup\{\epsi_i+\epsi_j,~1\leq i<j\leq n\}\cup\{2\epsi_i,~1\leq i\leq n\},\\
D_n^+&=\{\epsi_i-\epsi_j,~1\leq i<j\leq n\}\cup\{\epsi_i+\epsi_j,~1\leq i<j\leq n\}.\\
\end{split}
\end{equation*}
Here $\{\epsi_i\}_{i=1}^n$ is the standard basis of $\Rp^n$.

Denote by $\nt$ the algebra of all strictly upper-triangular matrices from $\gt$. Then $\nt$ has a basis consisting of root vectors $e_{\alpha}$, $\alpha\in\Phi^+$, where
\begin{equation*}\predisplaypenalty=0
\begin{split}
e_{\epsi_i}&=\sqrt{2}(e_{0,i}-e_{-i,0}),~e_{2\epsi_i}=e_{i,-i},\\
e_{\epsi_i-\epsi_j}&=\begin{cases}
e_{i,j}&\text{for }A_{n-1},\\
e_{i,j}-e_{-j,-i}&\text{for }B_n,~C_n\text{ and }D_n,
\end{cases}\\
e_{\epsi_i+\epsi_j}&=\begin{cases}
e_{i,-j}-e_{j,-i}&\text{for }B_n\text{ and }D_n,\\
e_{i,-j}+e_{j,-i}&\text{for }C_n,
\end{cases}
\end{split}
\end{equation*}
and $e_{i,j}$ are the usual elementary matrices. For $\ot_{2n}(\Cp)$ (respectively, for $\ot_{2n+1}(\Cp)$ and $\spt_{2n}(\Cp)$) we index the rows (from left to right) and the columns (from top to bottom) of matrices by the numbers $1,\ldots,n,-n,\ldots,-1$ (respectively, by the numbers $1,\ldots,n,0,-n,\ldots,-1$ and $1,\ldots,n,-n,\ldots,-1$). Note that $\gt=\htt\oplus\nt\oplus\nt_-$, where $\nt_-=\langle e_{-\alpha},~\alpha\in\Phi^+\rangle_{\Cp}$, and, by definition, $e_{-\alpha}=e_{\alpha}^T$. (The superscript~$T$ always stands for transposed.) The set $\{e_{\alpha},~\alpha\in\Phi\}$ can be extended to a unique Chevalley basis of $\gt$.

Let $G$ be one of the following classical Lie groups: $\SL_n(\Cp)$, $\SO_{2n}(\Cp)$, $\SO_{2n+1}(\Cp)$ or $\Sp_{2n}(\Cp)$. The group $\SO_{2n}(\Cp)$ (respectively, $\SO_{2n+1}$ and $\Sp_{2n}(\Cp)$) is realized as the subgroup of $\SL_{2n}(\Cp)$ (respectively, of $\SL_{2n+1}(\Cp)$ and $\SL_{2n}(\Cp)$) which preserves the form $\beta$. Let $H$ (respectively, $B$ and $N$) be the set of all diagonal (respectively, upper-triangular and upper-triangular with 1 on the diagonal) matrices from~$G$. Then $H$ is a maximal torus of $G$, $B$ is a Borel subgroup of $G$ containing $H$, $N$ is the unipotent radical of $B$, and $\gt$ (respectively, $\htt$, $\bt$ and $\nt$) is the Lie algebra of $G$ (respectively, of $H$, $B$ and $N$).

Denote by $U(\nt)$ the enveloping algebra of $\nt$, and by $S(\nt)$ the symmetric algebra of $\nt$. Then $\nt$ and $S(\nt)$ are $B$-modules as $B$ normalizes $N$. Denote by $Z(\nt)$ the center of $U(\nt)$. It is well-known that the restriction of the symmetrization map $$\sigma\colon S(\nt)\to U(\nt),~x^k\mapsto x^k,~x\in\nt,~k\in\Zp_{\geq0},$$ to the algebra $S(\nt)^N$ of $N$-invariants is an algebra isomorphism between $S(\nt)^N$ and $Z(\nt)$.

We next present a canonical set of generators of $Z(\nt)$ (or, equivalently, of $S(\nt)^N$), whose description goes back to J. Dixmier, A. Joseph and B. Kostant \cite{Dixmier1}, \cite{Joseph1}, \cite{Kostant1}, \cite{Kostant2}. Denote by $\Bu$ the following subset of $\Phi^+$:
\begin{equation*}
\Bu=\begin{cases}\bigcup\nolimits_{1\leq i\leq[n/2]}\{\epsi_i-\epsi_{n-i+1}\}&\text{for }A_{n-1},\\
\bigcup\nolimits_{1\leq i\leq n/2}\{\epsi_{2i-1}-\epsi_{2i},~\epsi_{2i-1}+\epsi_{2i}\}&\text{for }B_n,~n\text{ even},\\
\bigcup\nolimits_{1\leq i\leq[n/2]}\{\epsi_{2i-1}-\epsi_{2i},~\epsi_{2i-1}+\epsi_{2i}\}\cup\{\epsi_n\}&\text{for }B_n,~n\text{ odd},\\
\bigcup\nolimits_{1\leq i\leq n}\{2\epsi_i\}&\text{for }C_n,\\
\bigcup\nolimits_{1\leq i\leq[n/2]}\{\epsi_{2i-1}-\epsi_{2i},~\epsi_{2i-1}+\epsi_{2i}\}&\text{for }D_n.\\
\end{cases}
\end{equation*}
Note that $\Bu$ is a maximal strongly orthogonal subset of $\Phi^+$, i.e., $\Bu$ is maximal with the property that if $\alpha,\beta\in\Bu$ then neither $\alpha-\beta$ nor $\alpha+\beta$ belongs to $\Phi^+$. We call $\Bu$ the \emph{Kostant cascade} of orthogonal roots in~$\Phi^+$.

We can consider $\Zp\Phi$, the $\Zp$-linear span of $\Phi$, as a subgroup of the group $\Xt$ of rational multiplicative characters of $H$ by putting $\pm\epsi_i(h)=h_{i,i}^{\pm1}$, where $h_{i,j}$ is the $i$-th diagonal element of a matrix $h\in H$. Recall that a vector $\lambda\in\Rp^n$ is called a \emph{weight} of $H$ if $c(\alpha,\lambda)=2(\alpha,\lambda)/(\alpha,\alpha)$ is an integer for any $\alpha\in\Phi^+$, where $(\cdot,\cdot)$ is the standard inner product on $\Rp^n$. A weight $\lambda$ is called \emph{dominant} if $c(\alpha,\lambda)\geq0$ for all $\alpha\in\Phi^+$. An element $a$ of an $H$-module is called an \emph{$H$-weight vector}, if there exists $\nu\in\Xt$ such that $h\cdot a=\nu(h)a$ for all $h\in H$. By \cite[Theorems 6, 7]{Kostant2}, there exist unique (up to scalars) prime polynomials $\xi_1,~\ldots,~\xi_m\in S(\nt)^N$, $m=|\Bu|$, such that each $\xi_i$ is an $H$-weight polynomial of a dominant weight $\mu_i$ belonging to the $\Zp$-linear span $\Zp\Bu$ of $\Bu$. A remarkable fact is that
\begin{equation}
\xi_1,~\ldots,~\xi_m\text{ are algebraically independent generators of }S(\nt)^N,\label{formula:xi_i}
\end{equation}
so $S(\nt)^N$ and $Z(\nt)$ are polynomial rings. We call $\xi_i$ the \emph{$i$-th canonical generator} of $S(\nt)^N$.\newpage

Let $\nt^*$ be the dual space of $\nt$, and let $\{e_{\alpha}^*,~\alpha\in\Phi^+\}$ be the basis of $\nt^*$ dual to the basis $\{e_{\alpha},~\alpha\in\Phi^+\}$ of $\nt$. Put $R=\left\{t=\sum\nolimits_{\beta\in\Bu}t_{\beta}e_{\beta}^*,~t_{\beta}\in\Cp^{\times}\right\}$ and denote by $X$ the union of all $N$-orbits in $\nt^*$ of elements of $R$. In fact, $X$ is a single $B$-orbit in $\nt^*$, and the $N$-orbits of two distinct point of $R$ are disjoint. Kostant \cite[Theorems 1.1, 1.3]{Kostant3} proves that $X$ is a Zariski dense subset of $\nt^*$, and for $t\in R$, up to scalar multiplication,
\begin{equation}
\xi_i(t)=\prod_{\beta\in\Bu}t_{\beta}^{r_{\mu_i}(\beta)},~1\leq i\leq m,\text{ where }\label{formula:xi_of_t}
r_{\mu_i}(\beta)=\dfrac{(\mu_i,\beta)}{(\beta,\beta)}.
\end{equation}

The following representation-theoretic description of $\xi_i$ is given in \cite{Panov1}, however, it can be found in slightly different terms also in \cite{LipsmanWolf1} and \cite{Joseph-Faucaunt-Millet1}. Let $\Pi=\{\alpha_1,\ldots,\alpha_n\}$ be the set of simple roots \cite[Table I--IV]{Bourbaki1}, and let $\vpi_i$, $1\leq i\leq n$, be the $i$-th fundamental dominant weight of $\Phi$. We define positive integers $k_i$ for $1\leq i\leq m=|\Bu|$ by the following rule:
\begin{center}\begin{tabular}{|l|l|l|}
\hline
$\Phi=A_{n-1}$&$m=[n/2]$&$k_i=1$ for $1\leq i\leq m$\\
\hline
$\Phi=B_n$&$m=n$ for $B_n$&$k_i=1$ for odd $i$,\\
or $D_n$&$m=n$ for $D_n$ when $n$ even&$k_i=2$ for even $i<m$,\\
&$m=n-1$ for $D_n$ when $n$ odd&$k_m=1$ for even $m$\\
\hline
$\Phi=C_n$&$m=n$&$k_i=1$ for $1\leq i\leq m$\\
\hline
\end{tabular}\quad.
\end{center}

Let $W$ be the Weyl group of $\Phi$. Denote by $w_0$ the unique element of $W$ such that $w_0(\alpha)\in-\Phi^+$ for all $\alpha\in\Phi^+$. Furthermore, set $\vpi_i'=\dfrac{(1-w_0)\vpi_i}{k_i}$. Then the weights $\vpi_i'$'s have the following form:
\begin{center}\begin{tabular}{|l|l|}
\hline
$\Phi=A_{n-1}$&$\vpi_i'=2\epsi_1+\ldots+2\epsi_{i-1}+\epsi_i$ for $1\leq i\leq m$\\
\hline
$\Phi=B_n$&$\vpi_i'=\begin{cases}2\epsi_1+\ldots+2\epsi_i&\text{for odd $i<m$},\\
\epsi_1+\ldots+\epsi_i&\text{otherwise}
\end{cases}$\\
\hline
$\Phi=C_n$&$\vpi_i'=2\epsi_1+\ldots+2\epsi_i$ for $1\leq i\leq m$\\
\hline
$\Phi=D_n$&$\vpi_i'=\begin{cases}2\epsi_1+\ldots+2\epsi_i&\text{for odd $i<m-1$},\\
2\epsi_1+\ldots+2\epsi_{n-1}&\text{for $i=n-2=m-1$ when $n$ is odd},\\
\epsi_1+\ldots+\epsi_{n-1}-\epsi_n&\text{for $i=n-1=m-1$ when $n$ is even},\\
\epsi_1+\ldots+\epsi_i&\text{otherwise}
\end{cases}$\\
\hline
\end{tabular}\quad.
\end{center}

Now, let $V_i$ be the $i$-th fundamental representation of the group $G$ for $1\leq i\leq m$, and let $V_i^*$ be its dual representation. Fix highest-weight vectors $v_i$ and $l_i$ respectively of $V_i$ and $V_i^*$. Let $S_i$ be the regular function on $G$ defined by $S_i(g)=l_i(g\cdot v_i)$.
By $\exp\colon\gt\to G$ we denote the usual exponential map. Using the $\gt$-invariant form $\gt\ni x,y\mapsto\tr{xy}$ we identify the space $\nt_-$ of all lower-triangular matrices from $\gt$ with $\nt^*$. Then
\begin{equation}
e_{\alpha}^*=\begin{cases}\label{formula:Killing}
e_{\alpha}^T&\text{if }\Phi=A_{n-1},\text{ or }\Phi=C_n\text{ and }\alpha=2\epsi_i,\\
e_{\alpha}^T/4&\text{if }\Phi=B_n\text{ and }\alpha=\epsi_i,\\
e_{\alpha}^T/2&\text{otherwise}.
\end{cases}
\end{equation}
To each regular function $S$ on $G$ we assign the sequence of regular functions $S^j$ on $\nt_-$ (or, equivalently, of elements of $S(\nt)$) defined as the coefficients in the expansion $$S(\exp{tx})=t^k(S^0(x)+tS^1(x)+t^2S^2(x)+\ldots),~k\geq0,~x\in\nt_-.$$ In particular, regular functions $S_i^j$ on $G$ are now defined.

\mtheo{\textup{\cite[Theorem 2.12]{Panov1}}\\\indent\textup{i)} If $k_i=2$ for $1\leq i\leq m$\textup{,} then $S_i^0$ is a square in $S(\nt)$.\\\indent\textup{ii) }Set $K_i=S_i^0$ if $k_i=1$\textup{,} and $(K_i)^2=S_i^0$ if $k_i=2$. Then $K_i$ is a prime $H$-weight polynomial of weight $\vpi_i'$ contained in $S(\nt)^N$ for $1\leq i\leq m$.}\newpage

\corop{After suitable reordering of indices\textup{,} we have \textup{(}up to scalar mul\-ti\-pli\-cation\textup{)} $\mu_i=\vpi_i'$ and $\xi_i=K_i$ for $1\leq i\leq m$.}{It is easy to check that $\vpi_i'$ is a dominant weight in $\Zp\Bu$. Then our claim follows from Kostant's characterization of $\xi_1$, $\ldots$, $\xi_m$ as the unique (up to scalars) prime polynomials in $S(\nt)^N$ which are $H$-weight polynomials of dominant weights belonging to $\Zp\Bu$.}

We now fix explicit expressions for $\xi_i$ in all classical root systems.

i) $\Phi=A_{n-1}$. Let $V=\Cp^n$ be the natural representation of $G=\SL_n(\Cp)$. Then $V_i=\bigwedge^iV$ for all $i$ from 1 to $n-1$. If $e_1,\ldots,e_n$ is the standard basis of $V$, and $e_1^*,\ldots,e_n^*$ is the dual basis of $V^*$, then $v_i=e_1\wedge\ldots\wedge e_i$ is a highest-weight vector of $V_i$, and $l_i=e_{n-i+1}^*\wedge\ldots\wedge e_n^*$ is a highest-weight vector of~$V_i^*$. This implies that $S_i^0(x)$ is the lower left $(i\times i)$-minor of a matrix $x\in\nt_-$. Using (\ref{formula:Killing}) we can set
\begin{equation}
\xi_i=\begin{vmatrix}\label{formula:Delta_i_A_n}
e_{1,n-i+1}&\ldots&e_{1,n-1}&e_{1,n}\\
e_{2,n-i+1}&\ldots&e_{2,n-1}&e_{2,n}\\
\vdots&\iddots&\vdots&\vdots\\
e_{i,n-i+1}&\ldots&e_{i,n-1}&e_{i,n}\\
\end{vmatrix}=\begin{vmatrix}
e_{\epsi_1-\epsi_{n-i+1}}&\ldots&e_{\epsi_1-\epsi_{n-1}}&e_{\epsi_1-\epsi_n}\\
e_{\epsi_2-\epsi_{n-i+1}}&\ldots&e_{\epsi_2-\epsi_{n-1}}&e_{\epsi_2-\epsi_n}\\
\vdots&\iddots&\vdots&\vdots\\
e_{\epsi_i-\epsi_{n-i+1}}&\ldots&e_{\epsi_i-\epsi_{n-1}}&e_{\epsi_i-\epsi_n}\\
\end{vmatrix}.
\end{equation}

ii) $\Phi=C_n$. Let $V=\Cp^{2n}$ be the natural representation of $G=\Sp_{2n}(\Cp)$. Then $V_i=\bigwedge^iV$ for all $i$ from 1 to $n-1=m-1$. If $e_1,\ldots,e_n,e_{-n},\ldots,e_{-1}$ is the standard basis of~$V$, and $e_1^*,\ldots,e_n^*,e_{-n}^*,\ldots,e_{-1}^*$ is the dual basis of $V^*$, then $v_i=e_1\wedge\ldots\wedge e_i$ is a highest-weight vector of~$V_i$, and $l_i=e_{-i}^*\wedge\ldots\wedge e_{-1}^*$ is a highest-weight vector of $V_i^*$. Consequently, $S_i^0(x)$ is again the lower left $(i\times i)$-minor of a matrix $x\in\nt_-$, and for $1\leq i\leq m-1$ we can set
\begin{equation}
\xi_i=\begin{vmatrix}
e_{\epsi_1+\epsi_i}&\ldots&e_{\epsi_1+\epsi_3}&e_{\epsi_1+\epsi_2}&2e_{2\epsi_1}\\
e_{\epsi_2+\epsi_i}&\ldots&e_{\epsi_2+\epsi_3}&2e_{2\epsi_2}&e_{\epsi_1+\epsi_2}\\
e_{\epsi_3+\epsi_i}&\ldots&2e_{2\epsi_3}&e_{\epsi_2+\epsi_3}&e_{\epsi_1+\epsi_3}\\
\vdots&\iddots&\vdots&\vdots&\vdots\\
2e_{2\epsi_i}&\ldots&e_{\epsi_3+\epsi_i}&e_{\epsi_2+\epsi_i}&e_{\epsi_1+\epsi_i}\\
\end{vmatrix}.\label{formula:Delta_i_C_n}
\end{equation}

We claim that $\xi_m$ can also be defined via formula (\ref{formula:Delta_i_C_n}) for $i=m$. To verify this it is suffices to check that $\xi_m$ is proportional to the $m$-th canonical generator for $\spt_{2n+2}(\Cp)$. But the latter is obvious because this generator is a prime $N$-invariant $H$-weight polynomial of weight $\vpi_m'\in\Zp\Bu$.

iii) $\Phi=D_n$. Let $V=\Cp^{2n}$ be the natural representation of $G=\SO_{2n}(\Cp)$. Then $V_i=\bigwedge^iV$ for all $i$ from 1 to $n-2$. If $e_1,\ldots,e_n,e_{-n},\ldots,e_{-1}$ is the standard basis of~$V$, and $e_1^*,\ldots,e_n^*,e_{-n}^*,\ldots,e_{-1}^*$ is the dual basis of $V^*$, then $v_i=e_1\wedge\ldots\wedge e_i$ is a highest-weight vector of~$V_i$, and $l_i=e_{-i}^*\wedge\ldots\wedge e_{-1}^*$ is a highest-weight vector of $V_i^*$. First, if $i\leq n-2$ is even, then $k_i=2$, so $S_i^0=(K_i)^2$. It follows that $K_i(x)$, for $x\in\nt_-$, is the Pfaffian of the skew-symmetric matrix obtained from the lower left $i\times i$ submatrix of $x$ by reordering the columns. Therefore, for even $i\leq n-2$ we set
\begin{equation}\predisplaypenalty=0
\xi_i^2=\pm\begin{vmatrix}
e_{\epsi_1+\epsi_i}&\ldots&e_{\epsi_1+\epsi_3}&e_{\epsi_1+\epsi_2}&0\\
e_{\epsi_2+\epsi_i}&\ldots&e_{\epsi_2+\epsi_3}&0&-e_{\epsi_1+\epsi_2}\\
e_{\epsi_3+\epsi_i}&\ldots&0&-e_{\epsi_2+\epsi_3}&-e_{\epsi_1+\epsi_3}\\
\vdots&\iddots&\vdots&\vdots&\vdots\\
0&\ldots&-e_{\epsi_3+\epsi_i}&-e_{\epsi_2+\epsi_i}&-e_{\epsi_1+\epsi_i}\\
\end{vmatrix}.\label{formula:Delta_i_D_n_pf}
\end{equation}
Our normalization is such that the term $e_{\epsi_1+\epsi_2}e_{\epsi_3+\epsi_4}\ldots e_{\epsi_{i-1}+\epsi_i}$ enters $\xi_i$ with coefficient 1.

If $i\leq n-2$ is odd and $x=\sum_{\alpha\in\Phi^+}x_{\alpha}e_{\alpha}^T\in\nt_-$, then
\begin{equation}
\begin{split}
S_i(\exp{tx})&=\begin{vmatrix}
tx_{\epsi_1+\epsi_i}+O(t^2)&\ldots&tx_{\epsi_1+\epsi_2}+O(t^2)&2t^2a_1(x)+O(t^3)\\
tx_{\epsi_2+\epsi_i}+O(t^2)&\ldots&2t^2a_2(x)+O(t^3)&-tx_{\epsi_1+\epsi_2}+O(t^2)\\
\vdots&\iddots&\vdots&\vdots\\
2t^2a_i(x)+O(t^3)&\ldots&-tx_{\epsi_2+\epsi_i}+O(t^2)&-tx_{\epsi_1+\epsi_i}+O(t^2)\\
\end{vmatrix}\\
&=t^k(S_i^0(x)+tS_i^1(x)+t^2S_i^2(x)+\ldots)
\end{split}
\end{equation}
for some $k\geq0$, where $a_s(x)=\sum\nolimits_{j=s+1}^nx_{\epsi_s-\epsi_j}x_{\epsi_s+\epsi_j}$. Let $\Uu^i$ be the $(i\times i)$-matrix defined by $(\Uu^i)_{a,b}=-(\Uu^i)_{a,b}=e_{\epsi_a+\epsi_{i-b+1}}$ for $a<i-b+1$, $(\Uu^i)_{a,i-a+1}=0$, and $\Uu_s^i$ be the matrix obtained from $\Uu^i$ by deleting the $(i-s+1)$-th row and column. Then $k=i+1$ and we can set
\begin{equation}
\xi_i=\sum\nolimits_{s=1}^ia_s\det\Uu_s^i,
\label{formula:Delta_i_D_n_net}
\end{equation}
where $a_s=\sum\nolimits_{j=s+1}^ne_{\epsi_s-\epsi_j}e_{\epsi_s+\epsi_j}$.

Next, assume $i=m$ (in other words, $i=n-1$ and $n$ is odd, or $i=n$ and $n$ is even). As the $m$-th canonical generator for $\ot_{2n+2}(\Cp)$ is a prime $N$-invariant $H$-weight polynomial of weight $\vpi_m'\in\Zp\Bu$, it follows that $\xi_m$ can be also defined via formula~(\ref{formula:Delta_i_D_n_pf}) for $i=m$.

Finally, assume $n$ is even and $i=m-1=n-1$. In this case $\xi_{m-1}$ can be defined by
\begin{equation*}
\xi_{m-1}^2=\pm\begin{vmatrix}
e_{\epsi_1-\epsi_n}&e_{\epsi_1+\epsi_{n-1}}&\ldots&e_{\epsi_1+\epsi_2}&0\\
e_{\epsi_2-\epsi_n}&e_{\epsi_2+\epsi_{n-1}}&\ldots&0&-e_{\epsi_1+\epsi_2}\\
\vdots&\vdots&\iddots&\vdots&\vdots\\
e_{\epsi_{n-1}-\epsi_n}&0&\ldots&-e_{\epsi_2+\epsi_{n-1}}&-e_{\epsi_1+\epsi_{n-1}}\\
0&-e_{\epsi_{n-1}-\epsi_n}&\ldots&-e_{\epsi_2-\epsi_n}&-e_{\epsi_1-\epsi_n}\\
\end{vmatrix}
\end{equation*}
(our normalization is such that the term $e_{\epsi_{n-1}-\epsi_n}\xi_{m-3}$ enters $\xi_{m-1}$ with coefficient 1) as it is easy to check that this polynomial is a prime $N$-invariant $H$-weight polynomial of weight $\vpi_{m-1}'$.

iv) Case $\Phi=B_n$. Let $V=\Cp^{2n+1}$ be the natural representation of $G=\SO_{2n+1}(\Cp)$. Then $V_i=\bigwedge^iV$ for all $i$ from 1 to $n-1$. If $e_1,\ldots,e_n,e_0,e_{-n},\ldots,e_{-1}$ is the standard basis of~$V$, and $e_1^*,\ldots,e_n^*,e_0^*,e_{-n}^*,\ldots,e_{-1}^*$ is the dual basis of $V^*$, then $v_i=e_1\wedge\ldots\wedge e_i$ is a highest-weight vector of~$V_i$, and $l_i=e_{-i}^*\wedge\ldots\wedge e_{-1}^*$ is a highest-weight vector of $V_i^*$.

First, if $i\leq n-1$ is even, then $k_i=2$, so $S_i^0=(K_i)^2$. It follows that $\xi_i$ can be defined via formula~(\ref{formula:Delta_i_D_n_pf}). Second, if $i=n=m$ is even, then, arguing as above, we can define $\xi_m$ again by formula (\ref{formula:Delta_i_D_n_pf}) for $i=n$. Next, if $i\leq n-1$ is odd and $x=\sum_{\alpha\in\Phi^+}x_{\alpha}e_{\alpha}^T\in\nt_-$, then
\begin{equation}
\begin{split}
S_i(\exp{tx})&=\begin{vmatrix}
tx_{\epsi_1+\epsi_i}+O(t^2)&\ldots&tx_{\epsi_1+\epsi_2}+O(t^2)
&2t^2b_1(x)+O(t^3)\\
tx_{\epsi_2+\epsi_i}+O(t^2)&\ldots
&2t^2b_2(x)+O(t^3)&-tx_{\epsi_1+\epsi_2}+O(t^2)\\
\vdots&\iddots&\vdots&\vdots\\
2t^2b_i(x)+O(t^3)
&\ldots&-tx_{\epsi_2+\epsi_i}+O(t^2)&-tx_{\epsi_1+\epsi_i}+O(t^2)\\
\end{vmatrix}\\
&=t^k(S_i^0(x)+tS_i^1(x)+t^2S_i^2(x)+\ldots)
\end{split}
\end{equation}
for some $k\geq0$, where $b_s(x)=\sum\nolimits_{j=s+1}^nx_{\epsi_s-\epsi_j}x_{\epsi_s+\epsi_j}+x_{\epsi_s}^2$. It turns out that $k=i+1$ and we can define $\xi_i$ via formula (\ref{formula:Delta_i_D_n_net}) for $b_s=\sum\nolimits_{j=s+1}^ne_{\epsi_s-\epsi_j}e_{\epsi_s+\epsi_j}+e_{\epsi_s}^2/4$. Finally, assume $i=n=m$ is odd. Arguing as above, we see that $\xi_m$ can be defined by
\begin{equation*}
\xi_m^2=\pm\begin{vmatrix}
e_{\epsi_1}&e_{\epsi_1+\epsi_n}&\ldots&e_{\epsi_1+\epsi_2}&0\\
e_{\epsi_2}&e_{\epsi_2+\epsi_n}&\ldots&0&-e_{\epsi_1+\epsi_2}\\
\vdots&\vdots&\iddots&\vdots&\vdots\\
e_{\epsi_n}&0&\ldots&-e_{\epsi_2+\epsi_n}&-e_{\epsi_1+\epsi_n}\\
0&-e_{\epsi_n}&\ldots&-e_{\epsi_2}&-e_{\epsi_1}\\
\end{vmatrix}
\end{equation*}
(our normalization is such that the term $e_{\epsi_n}\xi_{m-2}$ enters $\xi_m$ with coefficient 1).

For $A_{n-1}$ and $C_n$, denote $\Delta_i=\sigma(\xi_i)$ for $1\leq i\leq m$. We call $\Delta_i$ the \emph{$i$-th canonical generator} of~$Z(\nt)$. Since all $e_{\alpha}$ involved in $\xi_i$ (i.e., $e_{\alpha}$ which appear in a term of $\xi_i$) commute, we conclude that  $\Delta_i$ is defined as an element of  $U(\nt)$ again by formulas (\ref{formula:Delta_i_A_n}), (\ref{formula:Delta_i_C_n}) for $A_{n-1}$, $C_n$ respectively. For $B_n$ and $D_n$, then denote $P_{i/2}=\sigma(\xi_i)$ when $i$ is even, and $\Du_{(i+1)/2}=\sigma(\xi_i)$ when $i$ is odd. If $i$ is even, all $e_{\alpha}$ involved in $\xi_i$ commute, so we can conclude that  $P_i$ is defined as an element of $U(\nt)$ again by formula (\ref{formula:Delta_i_D_n_pf}).\newpage

\sst{Infinite-dimensional case} Let \label{sst:centers_ifd} $\slt_{\infty}(\Cp)$, $\ot_{\infty}(\Cp)$, $\spt_{\infty}(\Cp)$ be the three simple complex finitary countable dimensional Lie algebras as classified by A.~Baranov \cite{Baranov1}. Each of them can be described as follows (see for example \cite{DimitrovPenkov2}). Consider an infinite chain of inclusions $$\gt_1\subset\gt_2\subset\ldots\subset\gt_n\subset\ldots$$ of simple Lie algebras, where $\rk{\gt_n}=n$ and all $\gt_n$ are of the same type $A$, $B$, $C$ or $D$. Then the union $\gt=\bigcup\gt_n$ is isomorphic to $\slt_{\infty}(\Cp)$, $\ot_{\infty}(\Cp)$ or $\spt_{\infty}(\Cp)$. It is always possible to choose nested Cartan subalgebras $\htt_n\subset\gt_n$, $\htt_n\subset\htt_{n+1}$, so that each root space of $\gt_n$ is mapped to exactly one root space of $\gt_{n+1}$. The union $\htt=\bigcup\htt_n$ acts semisimply on $\gt$, and is by definition a \emph{splitting Cartan subalgebra} of~$\gt$. We have a root decomposition $\gt=\htt\oplus\bigoplus_{\alpha\in\Phi}\gt^{\alpha}$ where $\Phi$ is \emph{the root system of} $\gt$ with respect to $\htt$ and $\gt^{\alpha}$ are the \emph{root spaces}. The root system $\Phi$ is simply the union of the root systems of $\gt_n$ and equals one of the following infinite root systems:
\begin{equation*}
\begin{split}
A_{\infty}&=\pm\{\epsi_i-\epsi_j,~i,j\in\Zp_{>0},~i<j\},\\
B_{\infty}&=\pm\{\epsi_i-\epsi_j,~i,j\in\Zp_{>0},~i<j\}\cup\pm\{\epsi_i+\epsi_j,~i,j\in\Zp_{>0},~i<j\}\cup\pm\{\epsi_i,~i\in\Zp_{>0}\},\\
C_{\infty}&=\pm\{\epsi_i-\epsi_j,~i,j\in\Zp_{>0},~i<j\}\cup\pm\{\epsi_i+\epsi_j,~i,j\in\Zp_{>0},~i<j\}\cup\pm\{2\epsi_i,~i\in\Zp_{>0}\},\\
D_{\infty}&=\pm\{\epsi_i-\epsi_j,~i,j\in\Zp_{>0},~i<j\}\cup\pm\{\epsi_i+\epsi_j,~i,j\in\Zp_{>0},~i<j\}.
\end{split}
\end{equation*}

A \emph{splitting Borel subalgebra} of $\gt$ is a subalgebra $\bt$ such that for every~$n$, $\bt_n=\bt\cap\gt_n$ is a Borel subalgebra of $\gt_n$. It is well-known that any splitting Borel subalgebra is conjugate via $\mathrm{Aut}\,\gt$ to a splitting Borel subalgebra containing $\htt$. Therefore, in what follows we restrict ourselves to considering only such Borel subalgebras $\bt$.

Recall \cite{DimitrovPenkov1} that a linear order on $\{0\}\cup\{\pm\epsi_i\}$ is $\Zp_2$-\emph{linear} if multiplication by $-1$ reverses the order. By \cite[Proposition 3]{DimitrovPenkov1}, there exists a bijection between splitting Borel subalgebras of $\gt$ containing $\htt$ and certain linearly ordered sets as follows.
\begin{equation*}\predisplaypenalty=0
\begin{split}
&\text{For }A_{\infty}\text{: linear orders on }\{\epsi_i\};\\
&\text{for }B_{\infty}\text{ and }C_{\infty}\text{: }\Zp_2\text{-linear orders on }\{0\}\cup\{\pm\epsi_i\};\\
&\text{for }D_{\infty}\text{: }\Zp_2\text{-linear orders on }\{0\}\cup\{\pm\epsi_i\}\text{ with the property that}\\
&\text{a minimal positive element (if it exists) belongs to }\Zp_{>0}.
\end{split}
\end{equation*}
In the sequel we denote these linear orders by $\prec$. To write down the above bijection, denote $\teta_i=\epsi_i$, if $i\succ0$, and $\teta_i=-\epsi_i$, if $\epsi_i\prec0$ (for $A_{\infty}$, $\teta_i=\epsi_i$ for all $i$). Then put $\bt=\htt\oplus\nt$, where $\nt=\bigoplus\limits_{\alpha\in\Phi^+}\gt^{\alpha}$ and
\begin{equation*}\predisplaypenalty=0
\begin{split}
A_{\infty}^+&=\{\teta_i-\teta_j,~i,j\in\Zp_{>0},~\teta_i\succ\teta_j\},\\
B_{\infty}^+&=\{\teta_i-\teta_j,~i,j\in\Zp_{>0},~\teta_i\succ\teta_j\}\cup\{\teta_i+\teta_j,~i,j\in\Zp_{>0},
~\teta_i\succ\teta_j\}\cup\{\teta_i,~i\in\Zp_{>0}\},\\
C_{\infty}^+&=\{\teta_i-\teta_j,~i,j\in\Zp_{>0},~\teta_i\succ\teta_j\}\cup\{\teta_i+\teta_j,~i,j\in\Zp_{>0},~
\teta_i\succ\teta_j\}\cup\{2\teta_i,~i\in\Zp_{>0}\},\\
D_{\infty}^+&=\{\teta_i-\teta_j,~i,j\in\Zp_{>0},~\teta_i\succ\teta_j\}\cup\{\teta_i+\teta_j,~i,j\in\Zp_{>0},~\teta_i\succ\teta_j\}.\\
\end{split}
\end{equation*}

Our goal in this subsection is to describe the center $Z(\nt)$ of the enveloping algebra $U(\nt)$. Fix $\nt$, i.e., fix an order $\prec$ as above. Define the subset $\Nu\subseteq\Zp_{>0}$ by setting $\Nu=\bigcup_{k\geq0}\Nu_k$, where $\Nu_0=\varnothing$ and $\Nu_k$ for $k\geq1$ is defined inductively in the following table.
\begin{center}
\begin{tabular}{|l|l|}
\hline
$\Phi$&$\Nu_k$\\
\hline\hline
$A_{\infty}$&$\Nu_{k-1}\cup\{i,j\}$ if there exists a maximal element $\teta_i$\\
&and a minimal element $\teta_j$ of $\{\teta_s,~s\in\Zp_{>0}\setminus\Nu_{k-1}\}$,\\
&$\Nu_{k-1}$ otherwise\\
\hline
$C_{\infty}$&$\Nu_{k-1}\cup\{i\}$ if there exists a maximal element $\teta_i$ of $\{\teta_s,~s\in\Zp_{>0}\setminus\Nu_{k-1}\}$,\\
&$\Nu_{k-1}$ otherwise\\
\hline
$B_{\infty}$,&$\Nu_{k-1}\cup\{i,j\}$ if there exists a maximal element $\teta_i$ of $\{\teta_s,~s\in\Zp_{>0}\setminus\Nu_{k-1}\}$\\
$D_{\infty}$&and a maximal element $\teta_j$ of $\{\teta_s,~s\in\Zp_{>0}\setminus\left(\Nu_{k-1}\cup\{i\}\right)\}$,\\
&$\Nu_{k-1}$ otherwise\\
\hline
\end{tabular}
\end{center}

\exam{i) Let $\Phi=A_{\infty}$. If $\epsi_1\succ\epsi_3\succ\ldots\succ\epsi_4\succ\epsi_2$, then $\Nu=\Zp_{>0}$. If $\epsi_1\succ\epsi_2\succ\epsi_3\succ\ldots$, then $\Nu=\varnothing$. ii) Let $\Phi\neq A_{\infty}$ and $\epsi_1\succ\epsi_2\succ\ldots\succ0\succ\ldots\succ-\epsi_2\succ-\epsi_1$. Then $\Nu=\Zp_{>0}$.}

Now we can define the (possibly infinite) Kostant cascade corresponding to $\nt$. Namely, to each $\Nu_k$ such that $\Nu_{k-1}\subsetneq\Nu_k$, we assign the root
\begin{equation*}
\beta_k=\begin{cases}\teta_i-\teta_j,&\text{if $\Phi=A_{\infty}$ and $\Nu_k\setminus\Nu_{k-1}=\{i,j\}$, $i\succ j$,}\\
\teta_i+\teta_j,&\text{if $\Phi=B_{\infty}$ or $D_{\infty}$ and $\Nu_k\setminus\Nu_{k-1}=\{i,j\}$, $i\succ j$,}\\
2\teta_i,&\text{if $\Phi=C_{\infty}$ and $\Nu_k\setminus\Nu_{k-1}=\{i\}$,}\\
\end{cases}
\end{equation*}
and put $\Bu=\{\beta_k,~k\geq1,~\Nu_{k-1}\subsetneq\Nu_k\}$. Note that $\Bu$ is a strongly orthogonal subset of $\Phi^+$; however it is not necessarily maximal with this property.

\defi{The subset $\Bu$ is called the \emph{Kostant cascade} corresponding to $\nt$.}
\exam{i) If $\Phi=A_{\infty}$ and $\epsi_1\succ\epsi_3\succ\ldots\succ\epsi_4\succ\epsi_2$, then $\Bu=\{\epsi_1-\epsi_2,~\epsi_3-\epsi_4$, $\epsi_5-\epsi_6,~\ldots\}.$ ii) If $\Phi\neq A_{\infty}$ and $\epsi_1\succ\epsi_2\succ\ldots\succ0\succ\ldots\succ-\epsi_2\succ-\epsi_1$, then
\begin{equation*}
\Bu=\begin{cases}\{\epsi_1+\epsi_2,~\epsi_3+\epsi_4,~\epsi_5+\epsi_6,~\ldots\}&\text{for $B_{\infty}$ and $D_{\infty}$},\\
\{2\epsi_1,~2\epsi_2,~2\epsi_3,~\ldots\}&\text{for $C_{\infty}$}.\\
\end{cases}
\end{equation*}}

To each finite non-empty subset $M\subset\Zp_{>0}$, one can assign a root subsystem $\Phi_M$ of $\Phi$ and a subalgebra~$\nt_M$ of $\nt$ by putting $\Phi_M=\Phi\cap\langle\epsi_i,~i\in M\rangle_{\Rp},~\nt_M=\bigoplus_{\alpha\in\Phi_M^+}\gt^{\alpha}$, $\Phi_M^+=\Phi_M\cap\Phi^+$. Then the subsystem $\Phi_M$ is isomorphic to the root system $\Phi_n$ of $\gt_n$ for $n=|M|$; we denote this isomorphism by $j_M\colon\Phi_n\to\Phi_M,~\epsi_i\mapsto\teta_{a_i}$, where $M=\{a_1,\ldots,a_n\}$, $\teta_{a_1}\succ\ldots\succ\teta_{a_n}$. Besides, $\nt_M$ is isomorphic as a Lie algebra to the maximal nilpotent subalgebra $\nt_n$ of $\gt_n$ considered in the previous subsection. Note also that $\nt=\ilm\nt_M$. Here, for $M\subseteq M'$, the monomorphism ${i_{M,M'}\colon\nt_M\hookrightarrow\nt_{M'}}$ is just the inclusion. Further, it is easy to see that there exist isomorphisms $\phi_M\colon\nt_n\to\nt_M,~M\subset\Zp_{>0},~n=|M|$, such that, for $M\subseteq M'$, $i_{M,M'}\circ\phi_M$ is just the restriction of $\phi_{M'}$ to $\nt_n\subset\nt_{n+1}$, and, for $\alpha\in\Phi_n^+$, $\phi_M(e_{\alpha})$ is a root vector corresponding to the root $j_M(\alpha)$; we denote it by $f_{j_M(\alpha)}$.

We are now ready to write down a set of generators of $Z(\nt)$. Namely, suppose that $q\geq1$ and $|\Bu|\geq q$. Let $M$ be a finite subset of $\Zp_{>0}$ such that $\Nu_q\subseteq M$. The isomorphism $\phi_M$ gives rise to an isomorphism $U(\nt_n)\to U(\nt_M)$, $n=|M|$. Slightly abusing notation, we denote the images of $\Delta_q$ and $P_q$ (as elements of $U(\nt_n)$ whenever defined) in $U(\nt_M)$ by the same letters. Then for $A_{\infty}$ (respectively, for~$C_{\infty}$), $\Delta_q\in U(\nt_M)$ is given by formula (\ref{formula:Delta_i_A_n}) (respectively, (\ref{formula:Delta_i_C_n})) for $i=q$ with $f_{j_M(\alpha)}$ instead of $e_{\alpha}$. Similarly, for $B_{\infty}$ and $D_{\infty}$, $P_q\in U(\nt_M)$ is given by formula (\ref{formula:Delta_i_D_n_pf}) for $i=2q$ with $f_{j_M(\alpha)}$ instead of~$e_{\alpha}$. It is important that $\Delta_q$, $P_q\in U(\nt_M)$ depend only on $q$ but not on $M$. Moreover, it is clear from the finite-dimensional theory that $\Delta_q$ (respectively, $P_q$) belong to the center of $U(\nt)$ for $A_{\infty}$ and $C_{\infty}$ (respectively, for $B_{\infty}$ and $D_{\infty}$).

Our first main result is as follows.

\theop{If \label{theo:center_ifd} $\Phi=A_{\infty}$\textup{,} $C_{\infty}$\textup{(}respectively\textup{,} $\Phi=B_{\infty}$\textup{,} $D_{\infty}$\textup{)}\textup{,} then $\Delta_q$ \textup{(}res\-pec\-tively\textup{,}~$P_q$\textup{)}\textup{,} $q\leq|\Bu|$\textup{,} generate $Z(\nt)$ as an algebra. In particular\textup{,} $Z(\nt)$ is a polynomial ring in $|\Bu|$ variables.}{Let $\Phi=A_{\infty}$, $C_{\infty}$, and $t\in Z(\nt)$. We need to prove that $t$ is a polynomial in $\Delta_q$, $q\leq|\Bu|$. Let $M$ be the minimal subset of $\Zp_{>0}$ for which $t\in U(\nt_M)$ (in particular, if $\Phi=A_{\infty}$, then $|M|$ is even), and $k$ be the maximal number such that $\Nu_k\subseteq M$ and $\Nu_{k-1}\subsetneq\Nu_k$. If $M=\Nu_k$, there is nothing to prove as the center of $U(\nt_M)$ is generated by $\Delta_1$, $\ldots$, $\Delta_k$. Therefore, assume that $\Nu_k\subsetneq M$. More precisely, let
\begin{equation*}
\Nu_k=\begin{cases}\{i_1,~\ldots,~i_k,~j_1,\ldots,~j_k\},&\text{if }\Phi=A_{\infty},\\
\{i_1,~\ldots,i_k\},&\text{if }\Phi=C_{\infty},
\end{cases}
\end{equation*}
$M\setminus\Nu_k=\{s_1,\ldots,s_r\}$, $n=|M|$, where $\teta_{i_1}\succ\ldots\succ\teta_{i_k}\succ\teta_{s_1}\succ\ldots\succ\teta_{s_r}\succ\teta_{j_k}\succ\ldots\succ\teta_{j_1}$.
In the rest of the proof we show that the assumption $\Nu_k\subsetneq M$ is contradictory.

Put $l=n/2=k+r/2$ (respectively, $l=n=k+r$) for $A_{\infty}$ (respectively, for $C_{\infty}$) and extend the set $\{\Delta_1,~\ldots,~\Delta_k\}$ to a set $\{\Delta_1,~\ldots,~\Delta_k,~\Delta_{k+1}',~\Delta_{k+2}',~\ldots,~\Delta_l'\}$ of generators of $Z(\nt_M)$ by letting $\Delta_i'$ be the image of the $i$-th canonical generator of $Z(\nt_n)$ under the isomorphism $\phi_M\colon U(\nt_n)\to U(\nt_M)$. For instance, let $k=2$, $r=4$ for $A_{\infty}$, or $k=2$, $r=2$ for $C_{\infty}$. Then $\Delta_4'$ has the following form.
\begin{center}
\begin{tabular}{|c|c|}
\hline
$\Phi=A_{\infty}$, $k=2$, $r=4$&$\Phi=C_{\infty}$, $k=2$, $r=2$\\
\hline\hline&\\
$\begin{vmatrix}
f_{\teta_{i_1}-\teta_{s_3}}&f_{\teta_{i_1}-\teta_{s_4}}&f_{\teta_{i_1}-\teta_{j_2}}&f_{\teta_{i_1}-\teta_{j_1}}\\
f_{\teta_{i_2}-\teta_{s_3}}&f_{\teta_{i_2}-\teta_{s_4}}&f_{\teta_{i_2}-\teta_{j_2}}&f_{\teta_{i_2}-\teta_{j_1}}\\
f_{\teta_{s_1}-\teta_{s_3}}&f_{\teta_{s_1}-\teta_{s_4}}&f_{\teta_{s_1}-\teta_{j_2}}&f_{\teta_{s_1}-\teta_{j_1}}\\
f_{\teta_{s_2}-\teta_{s_3}}&f_{\teta_{s_2}-\teta_{s_4}}&f_{\teta_{s_2}-\teta_{j_2}}&f_{\teta_{s_2}-\teta_{j_1}}\\
\end{vmatrix}$&$\begin{vmatrix}
f_{\teta_{i_1}+\teta_{s_2}}&f_{\teta_{i_1}+\teta_{s_1}}&f_{\teta_{i_1}+\teta_{i_2}}&2f_{2\teta_{i_1}}\\
f_{\teta_{i_2}+\teta_{s_2}}&f_{\teta_{i_2}+\teta_{s_1}}&2f_{2\teta_{i_2}}&f_{\teta_{i_1}+\teta_{i_2}}\\
f_{\teta_{s_1}+\teta_{s_2}}&2f_{2\teta_{s_1}}&f_{\teta_{i_2+\teta_{s_1}}}&f_{\teta_{i_1}+\teta_{s_1}}\\
2f_{2\teta_{s_2}}&f_{\teta_{s_1}+\teta_{s_2}}&f_{\teta_{i_2+\teta_{s_2}}}&f_{\teta_{i_1}+\teta_{s_2}}\\
\end{vmatrix}$\\&\\
\hline
\end{tabular}
\end{center}

Now we write $t=t'+t''$, where $t'$ belongs to the ideal of $Z(\nt_M)$ generated by $\Delta_{k+1}'$, $\ldots$, $\Delta_l'$, and $t''$ belongs to the subalgebra of $Z(\nt_M)$ generated by $\Delta_1$, $\ldots$, $\Delta_k$. Note that $t''\in Z(\nt)$, hence $t'\in Z(\nt)$. Moreover, $t'\neq0$ by the assumption that $\Nu_k\subsetneq M$.

The definition of $k$ implies that if $\Phi=C_{\infty}$ then $\teta_{s_1}$ is not a $\prec$-maximal element of the set $\{\teta_s,~s\in\Zp_{>0}\setminus\Nu_k\}$. Similarly, if $\Phi=A_{\infty}$ then at least one of the following holds: $\teta_{s_1}$ is not a \break$\prec$-maximal element of $\{\teta_s,~s\in\Zp_{>0}\setminus\Nu_k\}$, or $\teta_{s_r}$ is not a $\prec$-minimal element of $\{\teta_s,~s\in\Zp_{>0}\setminus\Nu_k\}$. In the sequel we assume that the former condition is satisfied (the case of the latter condition is similar). Then, in both cases $\Phi=A_{\infty}$ and $\Phi=C_{\infty}$, there exists $s_0\in\Zp_{>0}\setminus M$ such that $\teta_{i_k}\succ\teta_{s_0}\succ\teta_{s_1}$. Denote $M_0=M\cup\{s_0\}$, $\nt_0=\nt_{M_0}$ and $N_0=\exp{\nt_0}$. We have $t'\in Z(\nt_0)$. Put $\xi_i'=\sigma_0^{-1}(\Delta_i')$, $i=k+1,\ldots,l$, where $\sigma_0\colon S(\nt_0)\to U(\nt_0)$ is the symmetrization map. Note that $d'=\sigma_0^{-1}(t')\neq0$ as $t'\neq0$, and that $d'$ belongs to the ideal of $S(\nt_0)^{N_0}$ generated by $\xi_{k+1}'$, $\ldots$, $\xi_l'$.

Given $\alpha\in\Phi^+$, denote by $f_{\alpha}^*$ the linear form on $\nt$ (or on $\nt_M$, if $\alpha\in\Phi_M$) such that $f_{\alpha}^*(f_{\beta})=\delta_{\alpha,\beta}$ (the Kronecker delta). Furthermore, let $R_0$ be the subset of $\nt_0^*$ consisting of all elements of the form $\lambda=\sum_{\beta\in\Bu_0}\lambda_{\beta}f_{\beta}^*,~\lambda_{\beta}\in\Cp^{\times}$, where
\begin{equation*}
\Bu_0=\begin{cases}\{\teta_{i_1}-\teta_{j_1},\ldots,\teta_{i_k}-\teta_{j_k},\teta_{s_0}-\teta_{s_r}, \teta_{s_1}-\teta_{s_{r-1}},\ldots,\teta_{s_{r/2-1}}-\teta_{s_{r/2+1}}\}&\text{for }A_{\infty},\\
\{2\teta_{i_1},\ldots,2\teta_{i_k},2\teta_{s_0},2\teta_{s_1},\ldots,2\teta_{s_r}\}&\text{for }C_{\infty}.
\end{cases}
\end{equation*}
Denote by $X_0$ the union of the coadjoint $N_0$-orbits of all linear forms from $R_0$. As pointed out in Subsection~\ref{sst:centers_finite}, $X_0$ is a Zariski dense subset of $\nt_0^*$ by a result of Kostant.

If $\Phi=A_{\infty}$ and $\lambda\in R_0$, then $\lambda(f_{\alpha})=0$ for each $f_{\alpha}$ from the last row of $\xi_i'$ for $i=k+1,~\ldots,~l$. Therefore, $\xi_i'(\lambda)=0$ for all $\lambda\in R_0$. Thus, $d'(\lambda)=0$ for all $\lambda\in R_0$. However, since $d'$ is $N_0$-invariant, we obtain that the restriction of $d'$ to $X_0$ equals zero. As $X_0$ is Zariski dense in $\nt_0^*$, we conclude that $d'=0$, and consequently that $t'=0$. This contradiction completes the proof for $\Phi=A_{\infty}$.

Here is an illustration of the vanishing of $\xi_i'$ on $R_0$. Let $k=2$, $r=4$. On the picture below the boxes from $\Bu_0$ are marked by $\otimes$, and the boxes corresponding to the variables involved in $\xi_4'$ are grey:
{\Autonumfalse
\begin{equation*}\predisplaypenalty=0
\mymatrix{\lNote{$i_1$}\Note{$i_1$}\pho& \Note{$i_2$}\Lft{2pt}\Bot{2pt}\pho& \Note{$s_0$}\pho& \Note{$s_1$}\pho& \Note{$s_2$}\pho& \Note{$s_3$}\gray\pho& \Note{$s_4$}\gray\pho& \Note{$j_2$}\gray\pho& \Note{$j_1$}\gray\otimes\\
\lNote{$i_2$}\pho& \pho& \Lft{2pt}\Bot{2pt}\pho& \pho& \pho& \gray\pho& \gray\pho& \gray\otimes& \gray\pho\\
\lNote{$s_0$}\pho& \pho& \pho& \Lft{2pt}\Bot{2pt}\pho& \pho& \pho& \otimes& \pho& \pho\\
\lNote{$s_1$}\pho& \pho& \pho& \pho&\Lft{2pt}\Bot{2pt}\pho& \gray\otimes& \gray\pho& \gray\pho& \gray\pho\\
\lNote{$s_2$}\pho& \pho& \pho& \pho& \pho&\Lft{2pt}\Bot{2pt}\gray\pho& \gray\pho& \gray\pho& \gray\pho\\
\lNote{$s_3$}\pho& \pho& \pho& \pho& \pho& \pho& \Lft{2pt}\Bot{2pt}\pho& \pho& \pho\\
\lNote{$s_4$}\pho& \pho& \pho& \pho& \pho& \pho& \pho& \Lft{2pt}\Bot{2pt}\pho& \pho\\
\lNote{$j_2$}\pho& \pho& \pho& \pho& \pho& \pho& \pho& \pho& \Lft{2pt}\Bot{2pt}\pho\\
\lNote{$j_1$}\pho& \pho& \pho& \pho& \pho& \pho& \pho& \pho& \pho\\
}\qquad
\xi_4'=\begin{vmatrix}
f_{\teta_{i_1}-\teta_{s_3}}&f_{\teta_{i_1}-\teta_{s_4}}&f_{\teta_{i_1}-\teta_{j_2}}&f_{\teta_{i_1}-\teta_{j_1}}\\
f_{\teta_{i_2}-\teta_{s_3}}&f_{\teta_{2_1}-\teta_{s_4}}&f_{\teta_{i_2}-\teta_{j_2}}&f_{\teta_{2_1}-\teta_{j_1}}\\
f_{\teta_{s_1}-\teta_{s_3}}&f_{\teta_{s_1}-\teta_{s_4}}&f_{\teta_{s_1}-\teta_{j_2}}&f_{\teta_{s_1}-\teta_{j_1}}\\
f_{\teta_{s_2}-\teta_{s_3}}&f_{\teta_{s_2}-\teta_{s_4}}&f_{\teta_{s_2}-\teta_{j_2}}&f_{\teta_{s_2}-\teta_{j_1}}\\
\end{vmatrix}.\end{equation*}}

If $\Phi=C_{\infty}$, we write $t'=t_1+t_2$ where $t_2$ belongs to the subalgebra of $Z(\nt_0)$ generated by\break $\Delta_1$, $\ldots$, $\Delta_{k}$, and $t_1\neq0$ belongs to the ideal of $Z(\nt_0)$ generated by all other generators of $Z(\nt_0)$. Then $t_2\in Z(\nt)$, whence $t_1\in Z(\nt)$. Denote $d_1=\sigma_0^{-1}(t_1)\neq0$, $d_2=\sigma_0^{-1}(t_2)$. By definition, $d_2$ belongs to the ideal of the algebra $S(\nt_0)^{N_0}$ generated by the first $k$ canonical generators of $S(\nt_0)^{N_0}$. At the same time, $d'$ depends only on the $f_{\alpha}$'s from $\nt_M$. Since $d_1\neq0$, there exists $\lambda\in R_0$ such that $d_1(\lambda)\neq0$. On the other hand, set $\lambda_0=\lambda-\lambda_{2\teta_{s_0}}f_{2\teta_{s_0}}^*$. Then $d'(\lambda)=d'(\lambda_0)$ and $d_2(\lambda)=d_2(\lambda_0)$ by the definitions of $d'$ and $d_2$. Therefore $d_1(\lambda)=d_1(\lambda_0)$. However, $d_1(\lambda_0)=0$, because if $i>k$, then the value of the $i$-th canonical generator of $S(\nt_0)^{N_0}$ on $\lambda_0$ equals zero. This contradiction completes the proof for $\Phi=C_{\infty}$.

Assume now that $\Phi=D_n$ and $t$ is a central element of $U(\nt)$. We claim that $t$ belongs to the subalgebra of $U(\nt)$ generated by $P_q$ for $q\leq|\Bu|$. Let $M\subset\Zp_{>0}$ be a minimal finite set for which\break $t\in U(\nt_M)$ and $|M|=2n$ is even, and let $k$ be maximal such that $\Nu_k\subseteq M$ and $\Nu_{k-1}\subsetneq\Nu_k$. Denote $M\setminus\Nu_k=\{i_{k+1},j_{k+1},\ldots,i_n,j_n\}$, where $\teta_{i_{k+1}}\succ\teta_{j_{k+1}}\succ\ldots\succ\teta_{i_n}\succ\teta_{j_n}$ (if $n=k$, then $M=\Nu_k$). Set $N_M=\exp\nt_M$,
\begin{equation*}
\begin{split}
&\Bu_M=\{\teta_{i_1}-\teta_{j_1},~\teta_{i_1}+\teta_{j_1},~\ldots,~\teta_{i_n}-\teta_{j_n},~\teta_{i_n}+\teta_{j_n}\},\\
&R_M=\{t=\sum\nolimits_{\beta\in\Bu}t_{\beta}f_{\beta}^*\mid t_{\beta}\neq0\text{ for all }\beta\}.\\
\end{split}
\end{equation*}
Let $X_M$ be the union of all $N_M$-orbits of elements from $R_M$. Then $X_M$ is Zariski dense in $\nt_M^*$.

Denote the canonical generators of $S(\nt_M)^{N_M}$ by $\xi_1$, $\ldots$, $\xi_{2n}$. Then $Z(\nt_M)$ is generated as an algebra by $P_i=\sigma_M(\xi_{2i})$, $\Du_i=\sigma_M(\xi_{2i-1})$ for $1\leq i\leq n$, where $\sigma\colon S(\nt_M)\to U(\nt_M)$ is the symmetrization map. Further, set for simplicity $d_i=\xi_{2i-1}$, $1\leq i\leq n$. Using (\ref{formula:Delta_i_D_n_net}) one checks that for $t\in R_M$
\begin{equation}
\begin{split}
d_s(t)&=t_{\teta_{i_1}+\teta_{j_1}}^2\ldots t_{\teta_{i_{s-1}}+\teta_{j_{s-1}}}^2t_{\teta_{i_s}-\teta_{j_s}}t_{\teta_{i_s}+\teta_{j_s}},~1\leq s\leq n-1,\\
d_n(t)&=\label{formula:d_D_n}t_{\teta_{i_1}+\teta_{j_1}}\ldots t_{\teta_{i_{n-1}}+\teta_{j_{n-1}}}t_{\teta_{i_n}-\teta_{j_n}}.
\end{split}
\end{equation}

Assume $b$ is a positive integer such that there exists $a\in\Zp_{>0}\setminus M$ satisfying $\teta_{j_b}\succ\teta_a$. We can express~$t$ as $t=t'+t''$, where $t'$ belongs to the ideal of $Z(\nt_M)$ generated by $\Du_1$, $\ldots$, $\Du_b$, and $t''$ belongs to the subalgebra of $Z(\nt_M)$ generated by all remaining generators of $Z(\nt_M)$. Let $d=\sigma_M^{-1}(t)$, $d'=\sigma_M^{-1}(t')$, $d''=\sigma_M^{-1}(t'')$, so $d=d'+d''$. If $t'\neq0$, then $d'\neq0$. Let $c$ be the minimal among all numbers from~$1$ to~$b$ such that the variable $d_c$ appears in $d'$. Since $d'\neq0$, there exists $\lambda=\sum_{\beta\in\Bu_M}\lambda_{\beta}f_{\beta}^*\in R_M$ for which $d'(\lambda)\neq0$. Consider the set $Y=\lambda+\Cp f_{\teta_{i_c}-\teta_{j_c}}^*$. Obviously, $Y$ is a one-dimensional affine subspace of~$\nt_M^*$, and the restriction of $d'$ to $Y$ is a nonzero polynomial in one variable. Clearly, we can choose $\lambda$ so that this polynomial is of positive degree.

Now, put $M_a=M\cup\{a\}$, $N_a=\exp\nt_{M_a}$, $\mu=\lambda+f_{\teta_{i_c}+\teta_a}^*\in\nt_{M_a}^*$. Pick $s\in\Cp$ and put also $g=\exp{(sf_{\teta_{j_c}+\teta_a})}\in N_a$, $\mu'=g\cdot\mu$. Set $\Phi_a=\Phi_{M_a}$. One can easily check that $\mu'(f_{\alpha})=\mu(f_{\alpha})$ for all $\alpha\in\Phi_a^+\setminus\{\teta_{i_c}-\teta_{j_c}\}$, and $\mu'(f_{\teta_{i_c}-\teta_{j_c}})=\mu(f_{\teta_{i_c}-\teta_{j_c}})+s$.
Since $d\in S(\nt_{M_a})^{N_a}$, we obtain $d(\mu')=d(\mu)$ for all $s\in\Cp$. On the other hand, from the definition of $d''$ we see that $d''(\mu')=d''(\mu)$. Therefore, $d'(\mu')=d'(\mu)$. Define $\mu''$ as the restriction of $\mu'$ to~$\nt_M$. Then $\mu''$ belongs to $Y$, and $d'(\mu'')=d'(\mu')=d'(\mu)=d'(\lambda)$ as $d'\in S(\nt_M)$. Thus, the restriction of $d'$ to $Y$ is constant, a contradiction. We conclude that $d'=0$, and consequently that $t'=0$.

The above implies that it is sufficient to show that $M=\Nu_k$. Indeed, if $M=\Nu_k$, then $b$ can be chosen as $k$, and by the above $t$ is a polynomial in $P_i$, $1\leq i\leq k$. Assume, to the contrary, that $M\neq\Nu_k$, so $n>k$. Then there exists $s_0\in\Zp_{>0}\setminus M$ such that $\teta_{j_k}\succ\teta_{s_0}\succ\teta_{j_{k+1}}$. As we already know, $t$ belongs to the subalgebra of $Z(\nt_M)$ generated by $P_1$, $\ldots$, $P_k$, and by $P_i$, $\Du_i$ for $k+1\leq i\leq n$. We can express $t$ as $t=\wt t+t_0$, where $\wt t\neq0$ lies in the ideal of $Z(\nt_M)$ generated by $P_i$, $\Du_i$ for $k+1\leq i\leq n$, and $t_0$ lies in the subalgebra of $Z(\nt_M)$ generated by $P_1$, $\ldots$, $P_k$. If $M_0=M\cup\{s_0\}$ and $\nt_0=\nt_{M_0}$, then $\wt t\in Z(\nt_0)$.

Next, we write $\wt t$ as $\wt t=\wt t_1+\wt t_2$, where $\wt t_2$ lies in the subalgebra of $Z(\nt_0)$ generated by its first $2k$ canonical generators, and $\wt t_1$ lies in the ideal of $Z(\nt_0)$ generated by the remaining generators of $Z(\nt_0)$. If $\sigma_0\colon S(\nt_0)\to U(\nt_0)$ is the symmetrization map, we put $\wt d=\sigma_0^{-1}(\wt t),~\wt d_1=\sigma_0^{-1}(\wt t_1),~\wt d_2=\sigma_0^{-1}(\wt t_2)$. 
Let $R_0$ be the subset of $\nt_0^*$ consisting of all elements of the form $\lambda=\sum_{\beta\in\Bu_0}\lambda_{\beta}f_{\beta}^*,~\lambda_{\beta}\in\Cp^{\times}$, where $\Bu_0$ is the Kostant cascade for $\nt_0$. Note that the $(2k+1)$-th root in $\Bu_0$ equals $\teta_{s_0}+\teta_{i_{k+1}}$.

Denote by $X_0$ the union of the coadjoint $N_0$-orbits of all linear forms from $R_0$. Then $X_0$ is a Zariski dense subset of $\nt_0^*$. Assume $\wt d_1\neq0$. It is easy to check that there exists $\lambda\in R_0$ such that $\wt d_1(\lambda)\neq0$. Set $\lambda_0=\lambda-\lambda_{\teta_{s_0}+\teta_{i_{k+1}}}f_{\teta_{s_0}+\teta_{i_{k+1}}}^*$. Then $\wt d(\lambda)=\wt d(\lambda_0)$ and $\wt d_2(\lambda)=\wt d_2(\lambda_0)$ by definition of $d'$ and $d_2$. Hence $\wt d_1(\lambda)=\wt d_1(\lambda_0)=0$, a contradiction. Thus, $\wt d_1=0$, so $\wt d=\wt d_2$ belongs to the subalgebra of $Z(\nt_0)$ generated by its first $2k$ canonical generators. Finally, consider the affine space $Z=\sum_{r=1}^k(\Cp e_{\teta_{i_r}-\teta_{j_r}}+\Cp e_{\teta_{i_r}+\teta_{j_r}})=\langle e_{\beta},~\beta\in\Bu_M\cap\Bu_0\rangle_{\Cp}$. Let $\wt d_Z$ be the restriction of $\wt d$ to~$Z$. Since $\wt d=\wt d_2$, if $\wt d\neq0$ then $\wt d_Z$ is a nonzero polynomial of positive degree. But it follows from the definition of $\wt d$ that $\wt d_Z$ is zero (see (\ref{formula:d_D_n})). This shows that $M=\Nu_k$, and the proof for $\Phi=D_{\infty}$ is complete.

The proof for $\Phi=B_{\infty}$ is similar and we skip it.}

\sect{Centrally generated ideals of $U(\nt)$}
\sst{Finite-dimensional case} Throughout \label{sst:ideals_fd} this subsection $\gt$ and $\nt$ are as in Subsection~\ref{sst:centers_finite}. By definition, an ideal $J\subseteq U(\nt)$ is \emph{primitive} if $J$ is the annihilator of a simple $\nt$-module. Here we describe all primitive centrally generated ideals of $U(\nt)$, i.e., all primitive ideals $J$ generated (as ideals) by its intersection $J\cap Z(\nt)$ with the center $Z(\nt)$ of $U(\nt)$.

In the 1960s A.~Kirillov, B.~Kostant and J.-M. Souriau discovered that the orbits of the coadjoint action play a crucial role in the representation theory of $B$ and $N$ (see, e.g., \cite{Kirillov1}, \cite{Kirillov2}). The orbit method has a number of applications in the theory of integrable systems, symplectic geometry, etc. Work of J. Dixmier, M. Duflo, M. Vergne, O. Mathieu, N. Conze and R. Rentschler led to the result that the orbit method provides a nice description of primitive ideals of the universal enveloping algebra of a nilpotent Lie algebra (in particular, of $\nt$). Let us describe this in detail.

To any linear form $\lambda\in\nt^*$ one can assign a bilinear form $\beta_{\lambda}$ on $\nt$ by putting $\beta_{\lambda}(x,y)=\lambda([x,y])$. A subalgebra $\pt\subseteq\nt$ is a~\emph{polarization of $\nt$ at} $\lambda$ if it is a maximal $\beta_{\lambda}$-isotropic subspace. By \cite{Vergne}, such a subalgebra always exists. Let $\pt$ be a polarization of $\nt$ at $\lambda$, and $W$ be the one-dimensional representation of $\pt$ defined by $x\mapsto\lambda(x)$. Then the induced representation $V=U(\nt)\otimes_{U(\pt)}W$ of $\nt$ is irreducible. Hence, the annihilator $J(\lambda)=\Ann_{U(\nt)}{V}$ is a primitive two-sided ideal of~$U(\nt)$. It turns out that $J(\lambda)$ depends only on $\lambda$ and not on the choice of polarization. Further, $J(\lambda)=J(\mu)$ if and only if the coadjoint $N$-orbits of $\lambda$ and $\mu$ coincide. Finally, the \emph{Dixmier map} $$\Du\colon\nt^*\to\Prim U(\nt),~\lambda\mapsto J(\lambda),$$ induces a homeomorphism between $\nt^*/N$ and $\Prim U(\nt)$, where the latter set is endowed with the Jacobson topology. (See \cite{Dixmier2}, \cite{Dixmier3}, \cite{BorhoGabrielRentschler} for the details.)

In addition, it is well known that the following conditions on an ideal $J\subset U(\nt)$ are equivalent \cite[Proposition 4.7.4, Theorem 4.7.9]{Dixmier3}:
\begin{equation}
\begin{split}
&\text{i) $J$ is \label{formula:conditions_Prim_fd}primitive;}\\
&\text{ii) $J$ is maximal;}\\
&\text{iii) the center of $U(\nt)/J$ is trivial;}\\
&\text{iv) $U(\nt)/J$ is isomorphic to a Weyl algebra of finitely many variables.}
\end{split}
\end{equation}
Recall that the Weyl algebra $\Au_r$ of $r$ variables is the unital associative algebra with generators $p_i$,~$q_i$ for $1\leq i\leq r$, and relations $[p_i,q_i]=1$, $[p_i,q_j]=0$ for $i\neq j$, $[p_i,p_j]=[q_i,q_j]=0$ for all $i,~j$. Furthermore, in condition (\ref{formula:conditions_Prim_fd}) we have $U(\nt)/J\cong\Au_r$ where $r$ equals one half of the dimension of the coadjoint $N$-orbit of $\lambda$, given that $J=J(\lambda)$.

Recall the definition of the Kostant cascade $\Bu$ (Subsection~\ref{sst:centers_finite}) and set
\begin{equation*}
\Bu'=\begin{cases}\Bu&\text{for }\Phi=A_{n-1},~n\text{ odd},\\
\Bu\setminus\{\epsi_m-\epsi_{n-m+1}\}&\text{for }\Phi=A_{n-1},~n\text{ even},~m=n/2,\\
\bigcup_{1\leq i<n/2}\{\epsi_{2i-1}+\epsi_{2i+1}\}&\text{for }\Phi=B_n\text{ or }\Phi=D_n,\\
\Bu\setminus\{2\epsi_n\}&\text{for }\Phi=C_n.
\end{cases}
\end{equation*}
To a map $\xi\colon\Bu\to\Cp$ we assign the linear form $f_{\xi}=\sum_{\beta\in\Bu}\xi(\beta)e_{\beta}^*\in\nt^*$. We call a form $f_{\xi}$ a \emph{Kostant form} if $\xi(\beta)\neq0$ for any $\beta\in\Bu'$.


Let $V$ be a simple $\nt$-module and $J=\Ann_{U(\nt)}{V}$ be the corresponding primitive ideal of $U(\nt)$. By a version of Schur's Lemma \cite{Dixmier4}, each central element of $U(\nt)$ acts on $V$ as a scalar operator. For $A_{n-1}$ and $C_n$, let $c_k$ be the scalar corresponding to $\Delta_k$. For $B_n$ and $D_n$, let $c_k$ (respectively,~$\wt c_k$) be the scalar corresponding to~$P_k$ (respectively, to $\Du_k$). Note that these scalars do not depend on $V$ and are determined by $J$. Denote by $J_c$ the ideal of $U(\nt)$ generated by all $\Delta_k-c_k$ (respectively, by all $P_k-c_k$ and $\Du_k-\wt c_k$) for $A_{n-1}$ and $C_n$ (respectively, for $B_n$ and $D_n$). Clearly, $J_c\subseteq J$. Further, since $Z(\nt)$ is a polynomial ring and the center of $U(\nt)/J$ is trivial, $J$ is centrally generated if and only if $J=J_c$.  Put $m'=|\Bu'|$, $m=|\Bu|$.

Our second main result is as follows.

\mtheo{Suppose $\Phi$ is of type $A_{n-1}$ or $C_n$. The following conditions on a primitive ideal $J\subset U(\nt)$ are equivalent\textup{:}
\begin{equation*}
\begin{split}
&\text{\textup{i)} $J$ is centrally generated \textup{(}or\textup{,} equivalently\textup{,} $J=J_c$\textup{)}};\\
&\text{\textup{ii)} the scalars $c_1,\ldots,c_{m'}$ are nonzero};\\
&\text{\textup{iii)} $J=J(f_{\xi})$ for a Kostant form $f_{\xi}\in\nt^*$}.\\
\end{split}
\end{equation*}
If these\label{theo:ideal_fd} conditions are satisfied\textup{,} then the map $\xi$ is reconstructed by $J$\textup{:}
\begin{equation}
\xi(\beta)=(-1)^{k+1}\label{formula_xi_via_c_fd}c_k/c_{k-1},
\end{equation}
where $c_0=1$\textup{,} and $\beta=\epsi_k-\epsi_{n-k+1}$ for $\Phi=A_{n-1}$\textup{,} $\beta=2\epsi_k$ for $\Phi=C_n$.}


We expect this theorem to be true also for $B_n$ and $D_n$.

Before we prove Theorem~\ref{theo:ideal_fd} we prove few lemmas. We define the maps $\row\colon\Phi^+\to\Zp$ and $\col\colon\Phi^+\to\Zp$ by putting $\row(\epsi_i-\epsi_j)=\row(\epsi_i+\epsi_j)=\row(2\epsi_i)=i$, $\col(\epsi_i+\epsi_j)=\col(2\epsi_j)=-j$, $\col(\epsi_i-\epsi_j)=j$. Let $\Ro_i=\{\alpha\in\Phi^+\mid\row(\alpha)=i\}$. For $\alpha\in\Phi^+$, set
\begin{equation*}
\begin{split}
A(\alpha)&=\begin{cases}\bigcup\nolimits_{j+1\leq k\leq n-i+1}\{\epsi_j-\epsi_k\},&\text{if $\Phi=A_{n-1}$, $\alpha=\epsi_i-\epsi_j$, $j<n-i+1$},\\
\bigcup\nolimits_{n-j+1\leq k\leq i-1}\{\epsi_k-\epsi_i\},&\text{if $\Phi=A_{n-1}$, $\alpha=\epsi_i-\epsi_j$, $j>n-i+1$},\\
\bigcup\nolimits_{i\leq k\leq j-1}\{\epsi_k+\epsi_j\}\cup\Ro_j,&\text{if $\Phi=C_n$, $\alpha=\epsi_i-\epsi_j$},\\
\bigcup\nolimits_{i\leq k\leq j-1}\{\epsi_k-\epsi_j\},&\text{if $\Phi=C_n$, $\alpha=\epsi_i+\epsi_j$},\\
\end{cases}\\
\Bu(\alpha)&=\{\alpha\}\cup\{\beta\in\Bu\mid\row(\beta)<\row(\alpha)\},\\
R(\alpha)&=\{\row(\gamma),~\gamma\in\Bu(\alpha)\},~C(\alpha)=\{\col(\gamma),~\gamma\in\Bu(\alpha)\}.
\end{split}
\end{equation*}
Define a matrix $\Uu$ with entries from $U(\nt)$ by the following rule.

\begin{center}
\begin{tabular}{|l|l|l|}
\hline
$\Phi$&Size of $\Uu$&$\Uu$\\
\hline
$A_{n-1}$&$n\times n$&$\Uu_{i,j}=e_{\epsi_i-\epsi_j}$ for $1\leq i<j\leq n$,\\
&&$\Uu_{i,j}=0$ otherwise\\
\hline
$C_n$&$2n\times2n$&$\Uu_{i,j}=-\Uu_{-j,-i}=e_{\epsi_i-\epsi_j}$, $\Uu_{i,-j}=\Uu_{j, -i}=e_{\epsi_i+\epsi_j}$ for $1\leq i<j\leq n$,\\
&&$\Uu_{i,-i}=2e_{\epsi_2i}$, $1\leq i\leq n$, $\Uu_{i,j}=0$ otherwise\\
\hline
\end{tabular}
\end{center}

Denote by $\Delta_{\alpha}$ the element of $U(\nt)$, which equals the minor of $\Uu$ with rows $R(\alpha)$ and columns $C(\alpha)$. Note that the variables involved in each $\Delta_{\alpha}$ commute. For example, let $\Phi=A_{n-1}$, $n=8$, $\alpha=\epsi_3-\epsi_4$. On the picture below $\alpha$ is marked by $\bullet$, the roots from $\Bu$ are marked by $\otimes$'s, and the roots $\gamma$ such that $e_{\gamma}$ is involved in $\Delta_{\alpha}$ are grey:
\begin{equation*}\predisplaypenalty=0
\mymatrix{
\pho& \Lft{2pt}\Bot{2pt}\pho& \pho& \gray\pho& \pho& \pho& \gray\pho& \gray\otimes\\
\pho& \pho& \Lft{2pt}\Bot{2pt}\pho& \gray\pho& \pho& \pho& \gray\otimes& \gray\pho\\
\pho& \pho& \pho& \Lft{2pt}\Bot{2pt}\gray\bullet& \pho& \otimes& \gray\pho& \gray\pho\\
\pho& \pho& \pho& \pho& \Lft{2pt}\Bot{2pt}\otimes& \pho& \pho& \pho\\
\pho& \pho& \pho& \pho& \pho& \Lft{2pt}\Bot{2pt}\pho& \pho& \pho\\
\pho& \pho& \pho& \pho& \pho& \pho& \Lft{2pt}\Bot{2pt}\pho& \pho\\
\pho& \pho& \pho& \pho& \pho& \pho& \pho& \Lft{2pt}\Bot{2pt}\pho\\
\pho& \pho& \pho& \pho& \pho& \pho& \pho& \pho\\
}\qquad
\Delta_{\epsi_3-\epsi_4}=\begin{vmatrix}
e_{\epsi_1-\epsi_4}&e_{\epsi_1-\epsi_7}&e_{\epsi_1-\epsi_8}\\
e_{\epsi_2-\epsi_4}&e_{\epsi_2-\epsi_7}&e_{\epsi_2-\epsi_8}\\
e_{\epsi_3-\epsi_4}&e_{\epsi_3-\epsi_7}&e_{\epsi_3-\epsi_8}\\
\end{vmatrix}.\end{equation*}

\lemmp{Let $\alpha\in\Phi^+\setminus\Bu$. If $\gamma\notin A(\alpha)$ then $[\Delta_{\alpha},e_{\gamma}]=0$\label{lemm:comm1_A_n}. If $\gamma\in A(\alpha)$ then $[\Delta_{\alpha},e_{\gamma}]=\pm\Delta_{\alpha+\gamma}$. More precisely\textup{,} for $A_{n-1}$\textup{,} if $j<n-i+1$ then $[\Delta_{\epsi_i-\epsi_j},e_{\epsi_j-\epsi_s}]=\Delta_{\epsi_i-\epsi_s}$ for all $\epsi_j-\epsi_j\in A(\epsi_i-\epsi_j)$, and if $j>n-i+1$ then $[e_{\epsi_r-\epsi_i},\Delta_{\epsi_i-\epsi_j}]=\Delta_{\epsi_r-\epsi_j}$ for all $\epsi_r-\epsi_i\in A(\epsi_i-\epsi_j)$. For $C_n$\textup{,} if $\alpha=\epsi_i\pm\epsi_j$ then $[\Delta_{\alpha},e_{\gamma}]=\mp\Delta_{\alpha+\gamma}$ for all $\gamma\in A(\alpha)$.}{Let $\Phi=A_{n-1}$ (the proof for $C_n$ is similar). Suppose $\alpha=\epsi_i-\epsi_j$. Consider the case $j<n-i+1$ (the case $j>n-i+1$ can be considered similarly). If $e_{\epsi_p-\epsi_q}$ is involved in $\Delta_{\alpha}$ then $p\in\{1,2,\ldots,i\}=R(\alpha)$ and $q\in\{j,n-i+2,n-i+3,\ldots,n\}=C(\alpha)$. Assume that $e_{\epsi_r-\epsi_s}$ and $e_{\epsi_p-\epsi_q}$ do not commute. Then either $s\in R(\alpha)$ or $r\in C(\alpha)$. Consider these two cases separately.

First, if $s\in R(\alpha)$, i.e., $1\leq s\leq i$, then also $r\in R(\alpha)$ as $r<s$. Denote by $A_{p,q}$ the algebraic complement in $\Delta_{\alpha}$ to an element $e_{\epsi_p-\epsi_q}$. Then
$$\Delta_{\alpha}=e_{\epsi_s-\epsi_j}A_{s,j}+e_{\epsi_s-\epsi_{n-i+2}}A_{s,n-i+2}+e_{\epsi_s-\epsi_{n-i+3}}
A_{s,n-i+3}\ldots+e_{\epsi_s-\epsi_n}A_{s,n}.$$ Since $e_{\epsi_r-\epsi_s}$ commutes with each $e_{\gamma}$ involved in each algebraic complement, we have
\begin{equation*}
\begin{split}
[e_{\epsi_r-\epsi_s},\Delta_{\alpha}]&=[e_{\epsi_r-\epsi_s},e_{\epsi_s-\epsi_j}] A_{s,j}+[e_{\epsi_r-\epsi_s},e_{\epsi_s-\epsi_{n-i+2}}]
A_{s,n-i+2}+\ldots+[e_{\epsi_r-\epsi_s},e_{\epsi_s-\epsi_n}]A_{s,n}\\
&=e_{\epsi_r-\epsi_j}A_{s,j}+e_{\epsi_r-\epsi_{n-i+2}}A_{\epsi_s-\epsi_{n-i+2}}+\ldots+e_{\epsi_r-\epsi_n}A_{s,n}.\\
\end{split}
\end{equation*}
In other words, $[e_{\epsi_r-\epsi_s},\Delta_{\alpha}]$ equals the minor obtained from $\Delta_{\alpha}$ by replacing the $s$-th row by the $r$-th row. Thus $[e_{\epsi_r-\epsi_s},\Delta_{\alpha}]=0$.

Second, assume $r\in C(\alpha)=\{j,n-i+2,n-i+3,\ldots,n\}$. Clearly,
\begin{equation*}
\Delta_{\alpha}=e_{\epsi_1-\epsi_r}A_{1,r}+e_{\epsi_2-\epsi_r}A_{2,r}+\ldots+e_{\epsi_i-\epsi_r}A_{i,r},
\end{equation*}
so
\begin{equation*}
\begin{split}
[\Delta_{\alpha},e_{\epsi_r-\epsi_s}]&=[e_{\epsi_1-\epsi_r},e_{\epsi_r-\epsi_s}]A_{1,r}+[e_{\epsi_2-\epsi_r},e_{\epsi_r-\epsi_s}] A_{2,r}+\ldots+[e_{\epsi_i-\epsi_r},e_{\epsi_r-\epsi_s}]A_{i,r}\\
&=e_{\epsi_1-\epsi_s}A_{1,r}+e_{\epsi_2-\epsi_s}A_{2,r}+\ldots+e_{\epsi_i-\epsi_s}A_{i,r}.\\
\end{split}
\end{equation*}
Hence $[\Delta_{\alpha},e_{\epsi_r-\epsi_s}]$ equals the minor obtained from $\Delta_{\alpha}$ by replacing its $r$-th column by the column $e_{\epsi_1-\epsi_s},e_{\epsi_2-\epsi_s},\ldots,e_{\epsi_i-\epsi_s}$. If $r\neq j$, then the latter column is a column of $\Delta_{\alpha}$, so the commutator $[\Delta_{\alpha},e_{\epsi_r-\epsi_s}]$ is zero. If $r=j$ (and so $\epsi_r-\epsi_s=\epsi_j-\epsi_s\in A(\alpha)$), then the commutator equals $\Delta_{\epsi_i-\epsi_s}$ as required.}

Denote $\Delta_0=1$.

\lemmp{Let $\alpha,\beta\in\Phi^+$. Then $[\Delta_{\alpha},\Delta_{\beta}]=0$\textup{,} except\label{lemm:comm2_A_n} the following cases\textup{:}
\begin{equation*}
\begin{split}
[\Delta_{\epsi_i-\epsi_j},\Delta_{\epsi_j-\epsi_{n-i+1}}]&=(-1)^{i+1}\Delta_i\Delta_{i-1}\text{ if $\Phi=A_{n-1}$, $j<n-i+1$},\\
[\Delta_{\epsi_i-\epsi_j},\Delta_{\epsi_i+\epsi_j}]&=(-1)^{i+1}\Delta_i\Delta_{i-1}\text{ if $\Phi=C_n$}.\\
\end{split}
\end{equation*}\vspace{-0.5cm}}
{Suppose $\Phi=A_{n-1}$ (the proof for $C_n$ is similar). Consider the case $j<n-i+1$ (the case $j>n-i+1$ can be considered similarly). If $\beta=\epsi_r-\epsi_s\in\Bu$, then $\Delta_{\beta}=\Delta_r$ belongs to $Z(\nt)$, hence $\Delta_{\alpha}$ and $\Delta_{\beta}$ commute. So we may assume that $\beta\notin\Bu$, i.e., $s\neq n-r+1$. According to Lemma~\ref{lemm:comm1_A_n}, if $[\Delta_{\alpha},e_{\gamma}]\neq0$ then $\gamma\in A(\alpha)$. This implies that if $[\Delta_{\alpha},\Delta_{\beta}]\neq0$ then $\beta\in A(\alpha)$, because if $\beta\notin A(\alpha)$ then no $e_{\gamma}$ involved in $\Delta_{\beta}$ are contained in $A(\alpha)$. Hence $\beta=\epsi_j-\epsi_s$ for some $s$ such that $j+1\lee s\lee n-i+1$.

Suppose $\Delta_{\alpha}$ and $\Delta_{\beta}$ do not commute. Then, arguing as above, we see that $\alpha=\epsi_i-\epsi_j\in A(\beta)$. If $s<n-j+1$ then $B(\beta)$ consists of certain roots of the form $\epsi_s-\epsi_t$ for $s<t$, but $i<j<s$ so $\alpha\notin A(\beta)$. If $s=n-j+1$, then $\Delta_{\beta}=\Delta_j$ is a central element of $U(\nt)$, so it commutes with $\Delta_{\alpha}$. Finally, if $s>n-j+1$, then $A(\beta)$ consists of certain roots of the form $\epsi_k-\epsi_j$, $n-s+1\leq k\leq j-1$. Hence $n-s+1\leq i$, but $s\leq n-i+1$, so $s=n-i+1$. Thus, if $\beta\neq\epsi_j-\epsi_{n-i+1}$ then $[\Delta_{\alpha},\Delta_{\beta}]=0$.

It remains to compute $[\Delta_{\alpha},\Delta_{\epsi_j-\epsi_{n-i+1}}]$. One has
$$\Delta_{\alpha}=e_{\epsi_1-\epsi_j}A_{1,j}+e_{\epsi_2-\epsi_j}A_{2,j}+\ldots+e_{\epsi_i-\epsi_j}A_{i,j}.$$
The minor $\Delta_{\epsi_j-\epsi_{n-i+1}}$ commutes with all variables involved in this expression except for $e_{\epsi_i-\epsi_j}$. Since $A_{i,j}=(-1)^{i+1}\Delta_{i-1}$, we obtain
\begin{equation*}
\begin{split}
[\Delta_{\alpha},\Delta_{\epsi_j-\epsi_{n-i+1}}]=[(-1)^{i+1} e_{\epsi_i-\epsi_j}\Delta_{i-1},\Delta_{\epsi_j-\epsi_{n-i+1}}]=(-1)^{i+1}[e_{\epsi_i-\epsi_j},
\Delta_{\epsi_j-\epsi_{n-i+1}}]\Delta_{i-1}.
\end{split}
\end{equation*}
By Lemma~\ref{lemm:comm1_A_n}, $$[e_{\epsi_i-\epsi_j},\Delta_{\epsi_j-\epsi_{n-i+1}}]=\Delta_{\epsi_i-\epsi_{n-i+1}}=\Delta_i.$$ This concludes the proof.}

\lemmp{Suppose that $c_i\neq0$ for $1\leq i\leq m'=|\Bu'|$. Then the ideal~$J_c$\label{lemm:max_ideal_A_n} is primitive.}{Consider the case $\Phi=A_{n-1}$ (the proof for $C_n$ is similar). Put $\Au=U(\nt)/J_c$. Given $x\in U(\nt)$, denote by $\wt x$ its image in $\Au$ under the canonical projection. There is a natural partial order on $\Phi^+$: $\alpha>\beta$ if $\alpha-\beta$ is a sum of positive roots. Note that for $k\geq2$ we have
$$\Delta_k=\pm e_{\epsi_k-\epsi_{n-k+1}}\Delta_{k-1}+\text{terms containing only $e_{\alpha}$ for }\alpha>\epsi_k-\epsi_{n-k+1}.$$
It follows that in $\Au$ we can write any $\wt e_{\beta}$ for $\beta\in\Bu$ as a polynomial in $\wt e_{\alpha}$ for $\alpha\in\Phi^+\setminus\Bu$. In other words, $\Au$~is generated as an algebra by $\wt e_{\alpha}$ for $\alpha\in\Phi^+\setminus\Bu$.

Similarly, given $\alpha=\epsi_i-\epsi_j\in\Phi^+\setminus\Bu$, we have
\begin{equation*}
\begin{split}
\Delta_{\alpha}&=\pm e_{\alpha}\Delta_{k-1}+\text{terms containing only $e_{\gamma}$ for }\gamma>\alpha,\text{ where}\\
k&=\begin{cases}
i,&\text{if }j<n-i+1,\\
n-j+1,&\text{if }j>n-i+1.
\end{cases}
\end{split}
\end{equation*}
This implies that in $\Au$ one can write $\wt e_{\alpha}$ as a polynomial in $\wt\Delta_{\gamma}$ for $\gamma\in\Phi^+\setminus\Bu$. Thus, $\wt\Delta_{\gamma}$ for $\gamma\in\Phi^+\setminus\Bu$ generate $\Au$ as an algebra.

Now, given $\alpha=\epsi_i-\epsi_j\in\Phi^+\setminus\Bu$, $j<n-i+1$, let
\begin{equation*}
p_{\alpha}=\wt\Delta_{\alpha},~q_{\alpha}=(-1)^{i+1}c_i^{-1}c_{i-1}^{-1}\wt\Delta_{\epsi_j-\epsi_{n-i+1}}.
\end{equation*}
Lemma~\ref{lemm:comm2_A_n} shows that $[p_{\alpha},q_{\gamma}]=0$ for $\alpha\neq\gamma$, $[p_{\alpha},p_{\gamma}]=[q_{\alpha},q_{\gamma}]=0$ for all $\alpha$, $\gamma$, and $[p_{\alpha},q_{\alpha}]=1$. Hence $\Au$ is a quotient algebra of the Weyl algebra $\Au_N$ for $$N=(n-2)+(n-4)+\ldots=\#\{\epsi_i-\epsi_j\in\Phi^+\setminus\Bu\mid j<n-i+1\}.$$ But the Weyl algebra $\Au_N$ is simple, and $\Au\neq0$, so $\Au\cong\Au_N$. Thus $J_c$ is primitive (see (\ref{formula:conditions_Prim_fd})).}\newpage

\textsc{Proof of Theorem~\ref{theo:ideal_fd}}. $\mathrm{(ii)}\Longrightarrow\mathrm{(iii)}.$ Put $c_0=1$ and define $\xi(\beta)$, $\beta\in\Bu$, by formula (\ref{formula_xi_via_c_fd}). Denote $\pt=\langle e_{\epsi_i-\epsi_j},~j\leq n-i+1\rangle_{\Cp}$ for $A_{n-1}$, and $\pt=\langle e_{\alpha},~\col(\alpha)<0\rangle_{\Cp}$ for $C_n$. By \cite[Theorem 1.1]{Panov2} and \cite[Theorem 1.1]{Ignatyev1} the space $\pt$ is a polarization of $\nt$ at $f_{\xi}$. Let $V_{\xi}$ be the simple $\nt$-module corresponding to the linear form $f_{\xi}$ and the polarization $\pt$, i.e., $V_{\xi}=U(\nt)\otimes_{U(\pt)}W_{\xi}$, where $W_{\xi}$ is a one-dimensional representation of $\pt$ defined by $x\mapsto f_{\xi}(x)$. Then $\Delta_k$ acts on $V_{\xi}$ via the scalar $c_k$ for $1\leq k\leq m$. Consequently $J=J_c\subseteq J(f_{\xi})$, and so $J=J(f_{\xi})$.

$\mathrm{(iii)}\Longrightarrow\mathrm{(i)}.$ Let $\pt$, $V_{\xi}$ be as in the previous paragraph. Let $\beta_k$ be the $k$-th root from $\Bu$, i.e.,\break $\beta_k=\epsi_k-\epsi_{n-k+1}$ for $A_{n-1}$, and $\beta_k=2\epsi_k$ for $C_n$. Then $c_k\neq0$ for $1\leq k\leq m'$. Moreover, $\xi(\beta)=(-1)^{k+1}c_k/c_{k-1}$ for $1\leq k\leq m$, where $c_0=1$, because $\Delta_k$ acts on $V_{\xi}$ via the scalar $(-1)^{k+1}\xi(\beta_1)\ldots\xi(\beta_k)$. By Lemma~\ref{lemm:max_ideal_A_n}, the ideal $J_c$ is primitive. Since $J_c\subseteq J=J(f_{\xi})$, we have $J=J_c=J(f_{\xi})$, so $J$ is centrally generated.

$\mathrm{(i)}\Longrightarrow\mathrm{(ii)}.$ Assume, to the contrary, that some scalars $c_k$ equal zero. Suppose that $i_1\leq m'$ is the minimal number such that $c_{i_1}=0$. Now, define inductively two (finite) sequences $\{i_j\}$ and $\{k_j\}$ of positive integers by the following rule. If $i_j$ is already defined and there exists $k$ such that $i_j<k\leq m$ and $c_k\neq0$, then set $k_j$ to be the minimal among all such $k$. Similarly, if $k_j$ is already defined and there exists $i$ such $k_j<i\leq m'$ and $c_i=0$, then set $i_{j+1}$ to be the minimal among all such $i$.

To each $j$ for which $i_j$ exists we assign the root
\begin{equation*}
\gamma_j=\begin{cases}\epsi_{i_j}-\epsi_{n-i_j},&\text{if $\Phi=A_{n-1}$},\\
\epsi_{i_j}+\epsi_{i_j+2},&\text{if $\Phi=C_n$ and $i_j<m'=m-1=n-1$},\\
\epsi_{n-1}-\epsi_n,&\text{if $\Phi=C_n$ and $i_j=m'$}.
\end{cases}
\end{equation*}
To each $j$ such that both $i_j$ and $k_j$ exist we assign the set of roots
\begin{equation*}
\Gamma_j=\begin{cases}\{\epsi_{i_j}-\epsi_{n-k_j+1},~\epsi_{k_j}-\epsi_{n-i_j+1}\},&\text{if $\Phi=A_{n-1}$},\\
\{\epsi_{i_j}+\epsi_{k_j}\},&\text{if $\Phi=C_n$}.
\end{cases}
\end{equation*}
Denote the lengths of the sequences $\{i_j\}$, $\{k_j\}$ by $l_I$, $l_K$ respectively, and put
\begin{equation*}
X=\begin{cases}
\left(\Bu\setminus\left(\bigcup_{j=1}^r\{\beta_{i_j},~\beta_{k_j}\}\cup\bigcup_{j=i_{r+1}}^m\beta_j\right)\right)\cup
\bigcup_{j=1}^r\Gamma_j\cup\{\gamma_{r+1}\},&\text{if }l_I=r+1,~l_K=r,\\
\left(\Bu\setminus\left(\bigcup_{j=1}^r\{\beta_{i_j},~\beta_{k_j}\}\right)\right)\cup\bigcup_{j=1}^r\Gamma_j,&\text{if }l_I=l_K=r.\\
\end{cases}
\end{equation*}

Let $\vfi\colon X\to\Cp$ be a map. Put $\mu_{\vfi}=\sum_{\alpha\in X}\vfi(\alpha)e_{\alpha}^*$. By \cite[6.6.9 (c)]{Dixmier3}, $\Delta_k-c_k'\in J(\mu_{\vfi})$, $1\leq k\leq m$, where $c_k'=\xi_k(f)=\sigma^{-1}(\Delta_k)(f)$. By the definition of $\Delta_k$ there exist at least two distinct maps $\vfi_1$, $\vfi_2$ such that $c_k'=c_k$ for $1\leq k\leq m$. It follows from \cite[Theorem 1.4]{Panov2} and from the proof of \cite[Theo\-rem~3.1]{Ignatyev2} that the orbits of $\mu_{\vfi_1}$ and $\mu_{\vfi_1}$ are disjoint, so $J(\mu_{\vfi_1})\neq J(\mu_{\vfi_2})$. On the other hand, both $J(\mu_{\vfi_1})$ and $J(\mu_{\vfi_2})$ contain $J=J_c$, and this contradicts the maximality of $J$. The equivalence of (i), (ii), (iii) is now proved. The fact that the map $\xi$ is reconstructed by $J$ via formula (\ref{formula_xi_via_c_fd}) follows from the proof of the implication $\mathrm{(iii)}\Longrightarrow\mathrm{(i)}$.\hfill$\square$

Recall that $\lambda\in\nt^*$ is \emph{regular} if the $N$-orbit $\Omega_{\lambda}\subset\nt^*$ of $\lambda$ has maximal possible dimension. It follows from \cite[Theorem 2.3]{Kostant1} that all Kostant forms are regular. Moreover, for $\Phi=A_{n-1}$, a form $\lambda\in\nt^*$ is a Kostant form if and only if it is regular. Since it is known that an orbit $\Omega_{\lambda}$ contains at most one Kostant form, Theorem~\ref{theo:ideal_fd} that the Dixmier map establishes a bijection between Kostant forms and centrally generated ideals of $U(\nt)$ for $\Phi=A_{n-1}$, $C_n$.

For $\Phi=A_{n-1}$, Theorem~\ref{theo:ideal_fd} implies that a primitive ideal $J(\lambda)$ of $U(\nt)$ is centrally generated if and only if it is ``minimal'' in the sense that the orbit $\Omega_{\lambda}$ has maximal dimension. This reminds us of Duflo's famous theorem that if $\at$ is a semi-simple Lie algebra, then any minimal primitive ideal of $U(\at)$ is centrally generated \cite{Duflo1}. However, for $\Phi=C_n$, this analogy no longer holds as there exist regular forms $\lambda$ such that $J(\lambda)$ is not centrally generated (due to the fact that not every regular form $\lambda$ is a Kostant form).\newpage

\sst{Infinite-dimensional case} Throughout this subsection we use the notation from Sub\-sec\-tion~\ref{sst:centers_ifd}. We now restrict ourselves to the case $\Nu=\Zp_{>0}$. This means that, up to isomorphism, $\nt$ can be chosen to correspond to the linear order $\epsi_1\succ\epsi_3\succ\epsi_5\succ\ldots\succ\epsi_6\succ\epsi_4\succ\epsi_2$ for $A_{\infty}$ (respectively, to the linear order
$\epsi_1\succ\epsi_2\succ\epsi_3\succ\ldots\succ0\succ\ldots\succ-\epsi_2\succ-\epsi_1$ for all other root systems). In particular, $\teta_i=e_i$ for all $i\in\Zp_{>0}$, and $f_{\alpha}=e_{\alpha}$ for all $M\subset\Zp_{>0}$, $\alpha\in\Phi_M^+$. For $\alpha\in\Phi^+$, denote by $e_{\alpha}^*$ the linear form on $\nt$ such that $e_{\alpha}^*(e_{\beta})=\delta_{\alpha,\beta}$ (the Kronecker delta) for all $\beta\in\Phi^+$. In this subsection we describe all centrally generated ideals of $U(\nt)$ for $A_{\infty}$ and $C_{\infty}$.

For our choice of $\nt$, the Kostant cascade has the following form:
\begin{equation*}
\Bu=\begin{cases}
\{\epsi_1-\epsi_2,~\epsi_3-\epsi_4,~\ldots\}&\text{for }A_{\infty},\\[-2pt]
\{\epsi_1+\epsi_2,~\epsi_3+\epsi_4,~\ldots\}&\text{for }B_{\infty}\text{ and }D_{\infty},\\[-2pt]
\{2\epsi_1,~2\epsi_2,~\ldots\}&\text{for }C_{\infty}.\\
\end{cases}
\end{equation*}
The forms $f_{\xi}=\sum_{\beta\in\Bu}\xi(\beta)e_{\beta}^*\in\nt^*$ for all maps $\xi\colon\Bu\to\Cp^{\times}$ are by definition the \emph{Kostant forms} on $\nt$.

Our goal is to construct a partial Dixmier map, which attaches to each Kostant form a primitive ideal of $U(\nt)$. As in Subsection~\ref{sst:ideals_fd}, define the maps $\row\colon\Phi\to\Zp$ and $\col\colon\Phi\to\Zp$ by putting\break $\row(\epsi_i-\epsi_j)=\row(\epsi_i+\epsi_j)=\row(2\epsi_i)=i$, $\col(\epsi_i+\epsi_j)=\col(2\epsi_j)=-j$, $\col(\epsi_i-\epsi_j)=j$, and set $\Ro_k=\{\alpha\in\Phi^+\mid\row(\alpha)=k\}$. Put $\pt=\langle e_{\alpha},~\alpha\in\Phi^+\setminus\Mo\rangle_{\Cp}$, where
\begin{equation*}
\Mo=\begin{cases}
\{\epsi_i-\epsi_j,~\text{$i$ odd, $j$ even, $j<i$}\}&\text{for }A_{\infty},\\[-2pt]
\Ro_i,~\text{$i$ even}&\text{for $B_{\infty}$ and $D_{\infty}$},\\[-2pt]
\{\epsi_i-\epsi_j,~1\leq i<j\leq n\}&\text{for }C_{\infty}.\\
\end{cases}
\end{equation*}
Put also $\pt_n=\pt\cap\nt_n$, where $\nt_n=\nt_M$ for $M=\{1,\ldots,n\}$. Fix a Kostant form $f=f_{\xi}$. By \cite[Theo\-rem~1.1]{Ignatyev1}, $\pt_n$ is a polarization of $\nt_n$ at the linear form $f_n=f\mathbin{\mid}_{\nt_n}$. Thus, $\pt=\ilm\pt_n$ is a polarization of $\nt$ at~$f$. Moreover, denote
\begin{equation}
V_{\xi}=U(\nt)\otimes_{U(\pt)}W,~V_{\xi}^n=U(\nt_n)\label{formula:V_V_m_ifd}\otimes_{U(\pt_n)}W^n,
\end{equation}
where $W$ (respectively, $W^n$) is the one-dimensional representation of $\pt$ (respectively, of $\pt_n$) given by $x\mapsto f_{\xi}(x)$. The $\nt_n$-modules $V_n^{\xi}$ are simple and form a natural chain whose union is $V_{\xi}$. Hence, $V_{\xi}$ is a simple $\nt$-module. We denote its annihilator in $U(\nt)$ by $J(f_{\xi})$.

\nota{Let $\Po(f_{\xi})$ be the set of all polarizations $\at$ of $\nt$ at $f_{\xi}$ such that $\at_n=\at\cap\nt_n$ is a polarization of $\nt_n$ at $f_n$ for large enough $n$. Define $V_{\xi,\at}$ and $V_{\xi,\at}^n$ by formula (\ref{formula:V_V_m_ifd}) in which $\pt$ and $\pt_n$ are replaced by $\at\in P(f_\xi)$ and $\at_n$ respectively. Then $V_{\xi,\at}=\ilm{V_{\xi,\at}^n}$ is a simple $\nt$-module. As the annihilator of $V_{\xi,\at}^n$ in $U(\nt_n)$ does not depend on $\at_n$, we conclude that the annihilator of $V_{\xi,\at}$ does not depend on~$\at$. This shows that $J(f_{\xi})$ can be defined via any polarization $\at\in\Po(f_{\xi})$.}

\lemmp{For $A_{\infty}$ and\label{lemm:max_ideal_A_C_ifd} $C_{\infty}$\textup{,} the primitive ideal $J(f_{\xi})$ is generated by $\Delta_k-c_k$ for $k\geq1$.}
{It follows from the definition of $\pt$ that $\Delta_k$ acts on $V_{\xi}$ and on each $V_{\xi}^n$ for $n\geq k$ via the scalar $c_k$. Let $J_n$ be the annihilator of $V_{\xi}^n$ in $U(\nt_n)$, $n\geq1$. Theorem~\ref{theo:ideal_fd} implies that $J_n$ is generated by $\Delta_k-c_k$ for $1\leq k\leq m$, where $m=[n/2]$ for $A_{\infty}$ and $m=n$ for $C_{\infty}$. Hence $J_n=J(f_{\xi})\cap U(\nt_n)$. The result follows.\looseness=-1}

Now let $V$ be a simple $\nt$-module and $J$ be its annihilator. By \cite{Dixmier4}, for $A_{\infty}$ and $C_{\infty}$ (respectively, for $B_{\infty}$ and $D_{\infty}$), each $\Delta_k$ (respectively, $P_k$) acts on $V$ via some scalar $c_k$ for $k\geq1$.

Our third main result is as follows.

\theop{Let $\Phi=A_{\infty}$\textup{,} $C_{\infty}$\textup{,} and $\nt$ be as above. The following conditions on a primitive ideal $J\subset U(\nt)$ are equivalent:
\begin{equation*}\postdisplaypenalty=10000
\begin{split}
&\text{\textup{i)} $J$ is centrally generated};\\[-2pt]
&\text{\textup{ii)} all scalars $c_k$ are nonzero};\\[-2pt]
&\text{\textup{iii)} $J=J(f_{\xi})$ for a Kostant form $f_{\xi}$}.\\
\end{split}
\end{equation*}
If these\label{theo:ideal_ifd} conditions are satisfied\textup{,} then the scalars $c_k$ reconstruct $\xi(\beta)$ exactly as in Theorem~\textup{\ref{theo:ideal_fd}}.}{
$\mathrm{(ii)}\Longrightarrow\mathrm{(iii)}$. Define $\xi$ by formula (\ref{formula_xi_via_c_fd}). Then $J(f_{\xi})\subseteq J$ by Lemma~\ref{lemm:max_ideal_A_C_ifd}. On the other hand, consider $J_n=J(f_{\xi})\cap U(\nt_n)$ for $n\geq1$. The ideal $J_n$ of $U(\nt_n)$ contains $\Delta_k-c_k$ for $1\leq k\leq m$, where $m=[n/2]$ for $A_{\infty}$ and $m=n$ for $C_{\infty}$. Hence, according to Theorem~\ref{theo:ideal_fd}, $J_n$ is a maximal ideal of $U(\nt_n)$ contained in $J\cap U(\nt_n)$. Thus, $J_n=J\cap U(\nt_n)$ for $n\geq1$, i.e., $J=J(f_{\xi})$.

$\mathrm{(iii)}\Longrightarrow\mathrm{(i)}$. Follows from Lemma~\ref{lemm:max_ideal_A_C_ifd}.

$\mathrm{(i)}\Longrightarrow\mathrm{(ii)}$. Assume, to the contrary, that some scalars $c_k$ are zero. Let $i_1$ be the minimal number for which $c_{i_1}=0$. Define two (possibly, infinite) sequences $\{i_j\}$ and $\{k_j\}$ of positive integers inductively by the following rule. If $i_j$ is already defined and there exists $k>i_j$ such that $c_k\neq0$, then set $k_j$ to be the minimal number for which $c_{k_j}\neq0$. Similarly, if $k_j$ is already defined and there exists $i>k_j$ such that $c_i=0$, then set $i_{j+1}$ to be the minimal number such that $c_{i_{j+1}}=0$.

To each $j$ for which $i_j$ exists we assign the root
\begin{equation*}
\gamma_j=\begin{cases}
\epsi_{2i_j-1}-\epsi_{2i_j+2},&\text{if }\Phi=A_{\infty},\\
\epsi_{i_j}+\epsi_{i_j+1},&\text{ if}\Phi=C_{\infty}.
\end{cases}
\end{equation*}
To each $j$ such that both $i_j$ and $k_j$ exist we assign the set of roots
\begin{equation*}
\Gamma_j=\begin{cases}
\{\epsi_{2i_j-1}-\epsi_{2k_j},~\epsi_{2k_j-1}-\epsi_{2i_j}\},&\text{if }\Phi=A_{\infty},\\
\{\epsi_{i_j}+\epsi_{k_j}\},&\text{if }\Phi=C_{\infty}.
\end{cases}
\end{equation*}

Next, we define the subset $X\subset\Phi^+$ as in the proof of Theorem~\ref{theo:ideal_fd}. Namely, let $\beta_k$ be the $k$-th root from $\Bu$ (i.e., $\beta_k=\epsi_{2k-1}-\epsi_{2k}$ for $A_{\infty}$, and $\beta_k=2\epsi_k$ for $C_{\infty}$). Denote the lengths of the sequences $\{i_j\}$, $\{k_j\}$ by $l_I$, $l_K$ respectively, and put
\begin{equation*}\predisplaypenalty=0
X=\begin{cases}
\left(\Bu\setminus\left(\bigcup_{j\leq r}\{\beta_{i_j},~\beta_{k_j}\}\cup\bigcup_{j\geq i_{r+1}}\beta_j\right)\right)\cup\bigcup_{j\leq r}\Gamma_j\cup\{\gamma_{r+1}\},&\text{if }l_I=r+1,~l_K=r,\\
\left(\Bu\setminus\left(\bigcup_{j\leq r}\{\beta_{i_j},~\beta_{k_j}\}\right)\right)\cup\bigcup_{j\leq r}\Gamma_j,&\text{if }l_I=l_K=r,\\
\left(\Bu\setminus\left(\bigcup_{j\geq1}\{\beta_{i_j},~\beta_{k_j}\}\right)\right)\cup\bigcup_{j\geq1}\Gamma_j,&\text{if }l_I=l_K=\infty.\\
\end{cases}
\end{equation*}

Let $\mu_{\vfi}=\sum_{\alpha\in X}\vfi(\alpha)e_{\alpha}^*$, where $\vfi\colon X\to\Cp^{\times}$ is a map. To each $\alpha\in\Phi^+$ we assign the subset $S(\alpha)\subset\Phi^+$ as follows:
\begin{equation*}
S(\alpha)=\begin{cases}
\bigcup\nolimits_{\epsi_l\succ\epsi_i}\{\epsi_l-\epsi_j\},&\text{if }\alpha=\epsi_i-\epsi_j,\\
\bigcup\nolimits_{l=i+1}^{j-1}\{\epsi_l+\epsi_j\}\cup\Ro_j,&\text{if }\Phi=C_{\infty},~\alpha=2\epsi_i.\\
\end{cases}
\end{equation*}
We then set $\Mo=\bigcup_{\beta\in X}\Mo_{\beta}$ where $\Mo_{\beta}=\{\gamma\in S(\beta)\mid \gamma,\beta-\gamma\notin\bigcup\nolimits\Mo_{\alpha}\}$, the latter union being taken over all $\alpha\in X$ such that $\row(\alpha)<\row(\beta)$. Note that if $\beta\in X$, $\alpha,\gamma\in\Phi^+$, $\alpha\notin\Mo$ and $\alpha+\gamma=\beta$, then $\gamma\in\Mo$. This implies that $\mu_{\vfi}([x,y])=0$ for all $x, y\in\at=\langle e_{\alpha},~\alpha\in\Phi^+\setminus\Mo\rangle_{\Cp}$. Moreover, it is easy to see that $\at$ is a subalgebra of $\nt$, hence we can consider the $\nt$-module $V_{\vfi}=U(\nt)\otimes_{U(\at)}W_{\vfi}$, where $W_{\vfi}$ is the one-dimensional representation of $\at$ given by $x\mapsto\mu_{\vfi}(x)$. Let $J_{\vfi}$ be the annihilator of $V_{\vfi}$ in $U(\nt)$ (we do not assert that $J_{\vfi}$ is a primitive ideal as we do not discuss the irreducibility of $V_{\vfi}$). One can check that the map $\vfi$ can be chosen so that $J\subseteq J_{\vfi}$, so we assume in the rest of the proof that this condition is satisfied.

Let $\alpha$ be the unique root from $\Ro_{i_1}\cap X$. Explicitly,
\begin{equation*}
\alpha=\begin{cases}\epsi_{2i_1-1}-\epsi_{2i_1+2},&\text{if $\Phi=A_{\infty}$ and $c_k=0$ for all $k\geq i_1$},\\
\epsi_{2i_1-1}-\epsi_{2k_1},&\text{if $\Phi=A_{\infty}$ and $k_1>i_1$ is the minimal such that $c_{k_1}\neq0$},\\
\epsi_{i_1}+\epsi_{i_1+2},&\text{if $\Phi=C_{\infty}$ and $c_k=0$ for all $k\geq i_1$},\\
\epsi_{i_1}+\epsi_{k_1},&\text{if $\Phi=C_{\infty}$ and $k_1>i_1$ is the minimal such that $c_{k_1}\neq0$}.\\
\end{cases}
\end{equation*}

Given $\gamma\in\Phi^+$, let $M$ be a finite subset of $\Zp_{>0}$ such that $\gamma\in\Phi_M^+$, $n=|M|$. Recall the definition of $j_M$ and $\phi_M$ from Subsection~\ref{sst:centers_ifd}. Let $\Delta_{\gamma}$ be the image in $U(\nt_M)$ of $\Delta_{j_M^{-1}(\gamma)}\in U(\nt_n)$ under the isomorphism $\phi_M$. Note that $\Delta_{\gamma}$ depends only on $\gamma$ and not on $M$. We now show that $\Delta_{\alpha}-c\in J$ for some scalar $c$. We prove this for $\Phi=A_{\infty}$ (the proof for $\Phi=C_{\infty}$ is similar). First, suppose $\col(\alpha)=2i_1+2$. By Lemma~\ref{lemm:comm1_A_n}, $\Delta_{\alpha}$ commutes with all $e_{\gamma}$, except $\gamma=\beta_{i_1}-\alpha$, and in the latter case $[\Delta_{\alpha},e_{\gamma}]=\pm\Delta_{\alpha+\gamma}=\pm\Delta_{\beta_{i_1}}=\pm\Delta_{i_1}$. Hence the image of $\Delta_{\alpha}$ belongs to the center of the image of $U(\nt)$ in the algebra $\End_{\Cp}V$. By \cite{Dixmier4}, there exist $c\in\Cp$ such that $\Delta_{\alpha}-c\in J$.

Now we apply induction on $\col(\alpha)$ (the base is $\col(\alpha)=2i_1+2$). By Lemma~\ref{lemm:comm1_A_n}, for any $\alpha$, $\Delta_{\alpha}$ commutes with all $e_{\gamma}$ for $\gamma\notin A(\alpha)$. On the other hand, if $\gamma\in A(\alpha)$ then $[\Delta_{\alpha},e_{\gamma}]=\pm\Delta_{\alpha+\gamma}$. But $\gamma+\alpha\in\Ro_{i_1}$ and $2i_1\leq\col(\gamma+\alpha)<\col(\alpha)$. If $\col(\gamma+\alpha)=2i_1$, then $\Delta_{\gamma+\alpha}=\Delta_{i_1}\in J$. By the inductive assumption, if $\col(\gamma+\alpha)\geq2i_1+2$ then there exists $c'\in\Cp$ such that $\Delta_{\alpha+\gamma}-c'\in J\subseteq J_{\vfi}$. It follows from the definition of $V_{\vfi}$ that for $\delta\in\Ro_{i_1}$, where $\col(\delta)\geq2i_1$, $\Delta_{\delta}$ acts on $V_{\vfi}$ via the scalar $c_{\delta}=\pm\mu_{\vfi}(e_{\delta})\prod_{k<i_1}\vfi(\beta_k)$. In other words, $\Delta_{\delta}-c_{\delta}\in J_{\vfi}$. In particular, $\Delta_{\alpha+\gamma}-c_{\alpha+\gamma}\in J_{\vfi}$. Since $c_{\delta}$ in uniquely determined by $J_{\vfi}$, we conclude that $c'=c_{\alpha+\gamma}=0$ because $\mu_{\vfi}(e_{\gamma+\delta})=0$. Thus the image of $\Delta_{\alpha}$ belongs to the center of the image of $U(\nt)$ in the algebra $\End_{\Cp}V$. By \cite{Dixmier4} there exists $c\in\Cp$ such that $\Delta_{\alpha}-c\in J$.

Further, we see that $c=\pm\vfi(\alpha)\prod_{k<i_1}\vfi(\beta_k)$ because $J\subseteq J_{\vfi}$ and $\mu_{\vfi}(e_{\alpha})=\vfi(\alpha)$. Note also that $\Delta_k$ acts on $V_{\vfi}$ be the scalar $c_k'=(-1)^{k+1}\prod_{i\leq k}\vfi(\beta_i)$. Thus there exist at least two maps $\vfi_1$, $\vfi_2$ from $X$ to $\Cp^{\times}$ such that $\vfi_1(\alpha)\neq\vfi_2(\alpha)$ and $c_k'=c_k$ for all $k\in\Zp_{>0}$. This implies that both $J_{\vfi_1}$ and $J_{\vfi_2}$ contain $J$, which contradicts the uniqueness of $c$. The proof is complete.}

Denote by $\Au_{\infty}$ the Weyl algebra with countably many generators $p_i$, $q_i$ for $i\in\Zp_{>0}$, and relations
\begin{equation}
[p_i,q_i]=1,~[p_i,q_j]=0\text{ for }i\neq j,~[p_i,p_j]=[q_i,q_j]=0\text{ for\label{formula:Weyl_ifd} all }i,j.
\end{equation}
Note that the center of $\Au_{\infty}$ is trivial because the center of $\Au_r$ is trivial for any $r\geq1$. Similarly, $\Au_{\infty}$ is a simple algebra. We have the following corollary (cf. (\ref{formula:conditions_Prim_fd})).

\corop{Let $\Phi=A_{\infty}$\textup{,} $C_{\infty}$\label{coro:ifd_ideals}\textup{,} $\nt$ be as above, and $J$ be a primitive centrally generated ideal of~$U(\nt)$. Then
\begin{equation*}
\begin{split}
&\text{\textup{i)} $J$ is maximal\textup{;}}\\
&\text{\textup{ii)} the center of $U(\nt)/J$ is trivial\textup{;}}\\
&\text{\textup{iii)} $U(\nt)/J$ is isomorphic to the Weyl algebra $\Au_{\infty}$.}
\end{split}
\end{equation*}}
{(i) By Theorem~\ref{theo:ideal_ifd}, $J=J_{f_{\xi}}$ for some Kostant form $f_\xi$. It follows from the proof of Lemma~\ref{lemm:max_ideal_A_C_ifd} that $J\cap U(\nt_n)$ is maximal for all $n\geq1$. Hence $J$ is maximal.

(ii) This follows immediately from (iii).

(iii) One can construct a set of generators of $U(\nt)/J$ satisfying (\ref{formula:Weyl_ifd}) as in the proof of Lemma~\ref{lemm:max_ideal_A_n}. Since $\Au_{\infty}$ is simple, we have $U(\nt)/J\cong\Au_{\infty}$.}

We expect Theorem~\ref{theo:ideal_ifd} and Corollary~\ref{coro:ifd_ideals} to hold also for $B_{\infty}$ and $D_{\infty}$. Finally, we note that Theo\-rem~\ref{theo:ideal_ifd} establishes a one-to-one correspondence between centrally generated primitive ideals of $U(\nt)$ for $\Phi=A_{\infty}$, $C_{\infty}$ and Kostant forms.

\medskip\textsc{Mikhail Ignatyev: Samara State University, Ak. Pavlova 1, 443011 Samara, Russia}
\emph{E-mail address}: \texttt{mihail.ignatev@gmail.com}

\medskip\textsc{Ivan Penkov: Jacobs University Bremen, Campus Ring 1, 28759 Bremen, Germany}
\emph{E-mail address}: \texttt{i.penkov@jacobs-university.de}

\end{document}